\documentclass[a4paper,notitlepage,10pt]{article}%  
\usepackage{amssymb}
\usepackage{graphicx}
\usepackage{amsmath}
\usepackage{amsthm}
\usepackage{amsfonts}%
\usepackage{subfigure} 
 
\providecommand{\U}[1]{\protect\rule{.1in}{.1in}}
\theoremstyle{plain}
\newtheorem{theorem}{Theorem}[section]

\newtheorem{corollary}[theorem]{Corollary}

\newtheorem{lemma}[theorem]{Lemma}

\newtheorem{proposition}[theorem]{Proposition}

\theoremstyle{definition}

\newtheorem{remark}[theorem]{Remark}
\numberwithin{equation}{section}
\numberwithin{theorem}{section}

\let\footnote=\endnote

\ifx\pdfoutput\relax\let\pdfoutput=\undefined\fi
\newcount\msipdfoutput
\ifx\pdfoutput\undefined\else
\ifcase\pdfoutput\else
\ifx\paperwidth\undefined\else
\ifdim\paperheight=0pt\relax\else\pdfpageheight\paperheight\fi
\ifdim\paperwidth=0pt\relax\else\pdfpagewidth\paperwidth\fi
\fi\fi\fi
\begin{document}

\title{Positive solutions of an elliptic Neumann problem with a sublinear indefinite
nonlinearity }
\author{U. Kaufmann\thanks{FaMAF, Universidad Nacional de C\'{o}rdoba, (5000)
C\'{o}rdoba, Argentina. \textit{E-mail address: }kaufmann@mate.uncor.edu} , H.
Ramos Quoirin\thanks{Universidad de Santiago de Chile, Casilla 307, Correo 2,
Santiago, Chile. \textit{E-mail address: }humberto.ramos@usach.cl} , K.
Umezu\thanks{Department of Mathematics, Faculty of Education, Ibaraki
University, Mito 310-8512, Japan. \textit{E-mail address: }%
kenichiro.umezu.math@vc.ibaraki.ac.jp}
\and \noindent}
\maketitle

\begin{abstract}
Let $\Omega\subset\mathbb{R}^{N}$ ($N\geq1$) be a bounded and smooth domain
and $a:\Omega\rightarrow\mathbb{R}$ be a sign-changing weight satisfying
$\int_{\Omega}a<0$. We prove the existence of a \textit{positive} solution
$u_{q}$ for the problem
%\marginpar{modified}%
\[
\left\{
\begin{array}
[c]{lll}%
-\Delta u=a(x)u^{q} & \mathrm{in} & \Omega,\\
\frac{\partial u}{\partial\nu}=0 & \mathrm{on} & \partial\Omega,
\end{array}
\right.  \leqno{(P_{a,q})}
\]
if $q_{0}<q<1$, for some $q_{0}=q_{0}(a)>0$. In doing so, we improve the
existence result previously established in \cite{krqu}. In addition, we
provide the asymptotic behavior of $u_{q}$ as $q\rightarrow1^{-}$. When
$\Omega$ is a ball and $a$ is radial, we give some explicit conditions on $q$
and $a$ ensuring the existence of a positive solution of $(P_{a,q})$. We also
obtain some properties of the set of $q$'s such that $(P_{a,q})$ admits a
solution which is positive on $\overline{\Omega}$. Finally, we present some
results on nonnegative solutions having 
\textit{dead cores}. Our
approach combines bifurcation techniques, \textit{a priori }bounds and the
sub-supersolution method. Several methods and results apply as well to the Dirichlet counterpart of $(P_{a,q})$. 

\end{abstract}

%\maketitle
\bigskip

\section{Introduction}

Let $\Omega$ be a bounded and smooth domain of $\mathbb{R}^{N}$ with $N\geq1$,
and $0<q<1$. The purpose of this article is to discuss the existence of
\textit{positive} solutions for the problem
\[
\left\{
\begin{array}
[c]{lll}%
-\Delta u=a(x)u^{q} & \mathrm{in} & \Omega,\\
\frac{\partial u}{\partial\nu}=0 & \mathrm{on} & \partial\Omega,
\end{array}
\right.  \leqno{(P_{a,q})}
\]
where $\Delta$ is the usual Laplacian in $\mathbb{R}^{N}$, and $\nu$ is the
outward unit normal to $\partial\Omega$. Throughout this article, unless
otherwise stated, we assume that $r>N$ and $a\in L^{r}\left(  \Omega\right)  $
is such that
\[
a\text{ changes sign and }\int_{\Omega}a<0.\leqno{(H_0)}
\]
We set $a^{\pm}:=\max(\pm a,0)$. Note that the change of sign in $a$ 
% \marginpar{\textit{modified}} 
means that 
$|\text{supp }a^{\pm}|>0$, where $|A|$ stands for the Lebesgue measure of
$A\subset\mathbb{R}^{N}$. By a \textit{nonnegative} \textit{solution} of
$(P_{a,q})$ we mean a function $u\in W^{2,r}\left(  \Omega\right)  $ (and
hence $u\in\mathcal{C}^{1}(\overline{\Omega})$) that satisfies the equation
for the weak derivatives and the boundary condition in the usual sense, and
such that $u\geq0$ in $\Omega$. If, in addition, $u>0$ in $\Omega$, then we
call it a \textit{positive solution} of $(P_{a,q})$. In this case, we shall
also say that $(q,u)$ is a \textit{positive solution} of $(P_{a,q})$. Let us
denote by $\mathcal{P}^{\circ}$ the interior of the positive cone of
$\mathcal{C}^{1}(\overline{\Omega})$, i.e.,
\[
\mathcal{P}^{\circ}:=\left\{  u\in\mathcal{C}^{1}(\overline{\Omega}):u>0\text{
on }\overline{\Omega}\right\}  .
\]
We observe that a positive solution of $(P_{a,q})$ \textit{need not}
	%\marginpar{\textit{modified}}
belong to $\mathcal{P}^{\circ}$ (see e.g. Remark \ref{q0} (i) below).

Very few works have been devoted to $(P_{a,q})$, the first and main one 
	% \marginpar{\textit{modified}} 
being \cite{BPT}, where the following results were established (see 
Theorem 2.1 and Lemmas 2.1 and 3.1 therein):

\setcounter{theorem}{-1}

\begin{theorem}
[Bandle-Pozio-Tesei \cite{BPT}]\label{tbpt} Let $0<q<1$ and $a$ be a
sign-changing H\"older continuous function on $\overline{\Omega}$. Then, the
following three assertions hold:

\begin{enumerate}
\item If $(P_{a,q})$ has a positive solution then $\int_{\Omega}a<0$.

\item If $\int_{\Omega}a<0$, then $(P_{a,q})$ has at least one nontrivial
nonnegative solution.

\item $(P_{a,q})$ has at most one solution in $\mathcal{P}^{\circ}$.
\end{enumerate}
\end{theorem}

Let us mention that some of the above results were extended in \cite{alama} to
a problem that is a linear perturbation of $(P_{a,q})$. However, no
\textit{sufficient} conditions for the existence of \textit{positive}
solutions have been provided in \cite{alama,BPT}. Let us point out that, due
to the non-Lipschitzian character of $u^{q}$ at $u=0$ and the change of sign
in $a$, the strong maximum principle does not apply to $(P_{a,q})$. As a
consequence, one cannot derive the positivity of nontrivial nonnegative
solutions of $(P_{a,q})$.

To the best of our knowledge, the first existence result on positive solutions
of $(P_{a,q})$ has been proved in our recent work \cite{krqu}. We recall it
now. Denoting by $\Omega_{+}=\Omega_{+}(a)$ the largest open subset of
$\Omega$ where $a>0$ \textit{a.e.}, let us consider the following condition:%
\[
\Omega_{+}\text{ has finitely many connected components and }\left\vert
(\text{supp }a^{+})\setminus\Omega_{+}\right\vert =0.\leqno{(H_1)}
\]

Under $(H_{1})$, we 
	%\marginpar{\textit{modified}} 
showed that \textit{every}
nontrivial nonnegative solution of $(P_{a,q})$ belongs to $\mathcal{P}^{\circ
}$ if $q$ is close enough to $1$ (see \cite[Theorem 1.7]{krqu}). This
positivity result was proved via a continuity argument inspired by
\cite[Theorem 4.1]{J} (see also \cite{K}), which is based on the fact that the
strong maximum principle applies to $(P_{a,q})$ if $q=1$.
%\marginpar{modified}
As a consequence, assuming in addition $(H_{0})$, we deduced that, for $q$
close enough to $1$, $(P_{a,q})$ has a solution in $\mathcal{P}^{\circ}$,
which is the \textit{unique} nontrivial nonnegative solution of $(P_{a,q})$
(see \cite[Corollary 1.8]{krqu}). Let us mention that, in general, uniqueness
of nonnegative solutions for $(P_{a,q})$ \textit{does not} hold (see e.g. the
proof of Theorem \ref{ti} (ii) below).

Regarding the Dirichlet counterpart of $(P_{a,q})$, we refer to
\cite{bandle,PT} for the existence of nontrivial nonnegative solutions, and to
\cite{nodea,ans,jesusultimo,krqu} for the existence of a positive solution. Let us mention, as already pointed out in \cite{alama,bandle,BPT}, that problems like $(P_{a,q})$ and its Dirichlet counterpart naturally arise in population dynamics models, cf. \cite{GM,N}. 

Our purpose in this article is to carry on the investigation of $(P_{a,q})$,
refining and extending % \marginpar{\textit{modified}} 
the existence results on positive solutions established in
\cite{krqu}.
%\marginpar{\textit{modified}}
In particular, following a different approach to 
 the one in \cite{krqu}, we 
shall remove $(H_{1})$ and prove that under $(H_{0})$ the problem $(P_{a,q})$
has a unique solution $u_{q}\in\mathcal{P}^{\circ}$ for $q$ close to $1$. As a byproduct, we deduce that $(H_0)$ is necessary and sufficient for the existence of a positive solution of $(P_{a,q})$ for {\it some} $q \in (0,1)$, see Corollary \ref{cor}. 
	% \marginpar{modified}
Moreover, we shall provide the stability properties of $u_{q}$ and its
asymptotic behavior as $q\rightarrow1^{-}$ (see Theorem \ref{c1} below). Note
that the stability analysis for solutions in $\mathcal{P}^{\circ}$ of
$(P_{a,q})$ is not easily carried out for $q\in(0,1)$ in general (see Remark
\ref{rem:stab} (ii)). 

Under $(H_0)$, 
let us denote by $\mu_{1}\left(  a\right)  $ the first positive eigenvalue of
the problem
\[
\left\{
\begin{array}
[c]{lll}%
-\Delta\phi=\mu a(x)\phi & \mathrm{in} & \Omega,\\
\frac{\partial\phi}{\partial\nu}=0 & \mathrm{on} & \partial\Omega,
\end{array}
\right.  \leqno{(E_{\mu,a})}
\]
and by $\phi_{1}=\phi_{1}(a)$ the associated positive eigenfunction satisfying
$\int_{\Omega}\phi_{1}^{2}=1$. It is well known that $\mu_{1}(a)$ is simple,
and $\phi_{1}\in\mathcal{P}^{\circ}$.

We shall look at $q$
as a bifurcation parameter in $(P_{a,q})$. As a matter of fact, note that if
$\mu_{1}(a)=1$, then $u=t\phi_{1}$ solves $(P_{a,1})$, i.e., $(P_{a,q})$ has
the trivial line $\Gamma_{1}$ of solutions in $\mathcal{P}^{\circ}$, where
\[
\Gamma_{1}:=\{(q,u)=(1,t\phi_{1}):t>0\}.
\]
We shall obtain, for $q$ close to $1$, 
	% \marginpar{modified}
a curve of solutions in $\mathcal{P}^{\circ}$ bifurcating from $\Gamma_{1}$ (see Figure \ref{figbif}).

Let us recall that a solution $u\in\mathcal{P}^{\circ}$ of $(P_{a,q})$ is said
to be \textit{asymptotically stable} (respect.\ \textit{unstable}) if
$\gamma_{1}(q,u)>0$ (respect.\ $<0$), where $\gamma_{1}(q,u)$ is the first
eigenvalue of the linearized eigenvalue problem at $u$, namely,%
\[
\left\{
\begin{array}
[c]{lll}%
-\Delta\varphi=qa(x)u^{q-1}\varphi+\gamma\varphi & \mathrm{in} & \Omega,\\
\frac{\partial\varphi}{\partial\nu}=0 & \mathrm{on} & \partial\Omega.
\end{array}
\right.
\]
In addition, $u$ is said to be \textit{weakly stable} if $\gamma_{1}%
(q,u)\geq0$.

Set
\begin{equation}
t_{\ast}:=\exp\left[  -\frac{\int_{\Omega}a\left(  x\right)  \phi_{1}^{2}%
\log\phi_{1}}{\int_{\Omega}a\left(  x\right)  \phi_{1}^{2}}\right]  .
\label{et*}%
\end{equation}
We are now in position to state our main results.

\begin{theorem}
\label{c1} Assume $(H_{0})$. Then there exists $q_{0}=q_{0}\left(  a\right)
\in(0,1)$ such that $(P_{a,q})$ has a unique solution $u_{q}\in\mathcal{P}%
^{\circ}$ for $q_{0}<q<1$. Moreover, $u_{q}$ is asymptotically stable and
satisfies the asymptotics
\[
u_{q} \sim\mu_{1}(a)^{-\frac{1}{1-q}} \, t_{\ast} \phi_{1} \quad\mbox{as}
\quad q\rightarrow1^{-},
\]
i.e. $\mu_{1}(a)^{\frac{1}{1-q}} u_{q} \rightarrow t_{*} \phi_{1}$ in
$\mathcal{C}^{1}(\overline{\Omega})$ as $q \to1^{-}$. More specifically:

\begin{enumerate}
\item If $\mu_{1}(a)=1$, then $u_{q}\rightarrow t_{\ast}\phi_{1}$ in
$\mathcal{C}^{1}(\overline{\Omega})$ as $q\rightarrow1^{-}$.

\item If $\mu_{1}(a)>1$, then $u_{q}\rightarrow0$ in $\mathcal{C}%
^{1}(\overline{\Omega})$ as $q\rightarrow1^{-}$.

\item If $\mu_{1}(a)<1$, then $\displaystyle\min_{\overline{\Omega}}%
u_{q}\rightarrow\infty$ as $q\rightarrow1^{-}$.
\end{enumerate}
\end{theorem}

One may easily see that Theorem \ref{tbpt} (i) still holds if $a \in
L^{r}(\Omega)$, with $r>N$, cf. the proof of \cite[Lemma 2.1]{BPT}. As a
consequence of this result and Theorem \ref{c1}, we derive the following:

\begin{corollary} \label{cor}
$(P_{a,q})$ has a positive solution (or a solution in $\mathcal{P}^{\circ}$)
for \textit{some} $q \in(0,1)$ if 
	% \marginpar{\textit{modified}} 
and only if $(H_{0})$ holds.
\end{corollary}

We shall prove Theorem \ref{c1} using a bifurcation technique based on the
Lyapunov-Schmidt reduction, which yields the existence of bifurcating
solutions in $\mathcal{P}^{\circ}$ from $\Gamma_{1}$ provided that $\mu
_{1}(a)=1$. By a suitable rescaling, we deduce then the results for the case
$\mu_{1}(a)\neq1$.
%\marginpar{\textit{added}}
Let us also 
	% \marginpar{\textit{modified}} 
point out that, in general, it is hard to give a lower estimate for
$q_{0}\left(  a\right)  $, see Remark \ref{q0} (i) below.

%%%%%%%%%%%%%%%%%%%%%%%%%%
\def\normu{\| u \|_{\mathcal{C}^1(\overline{\Omega})}}
\def\minu{\displaystyle{\min_{\overline{\Omega}}u}}
\def\bpt{(1, t_* \phi_{1})}
%
%	\begin{figure}[H]
	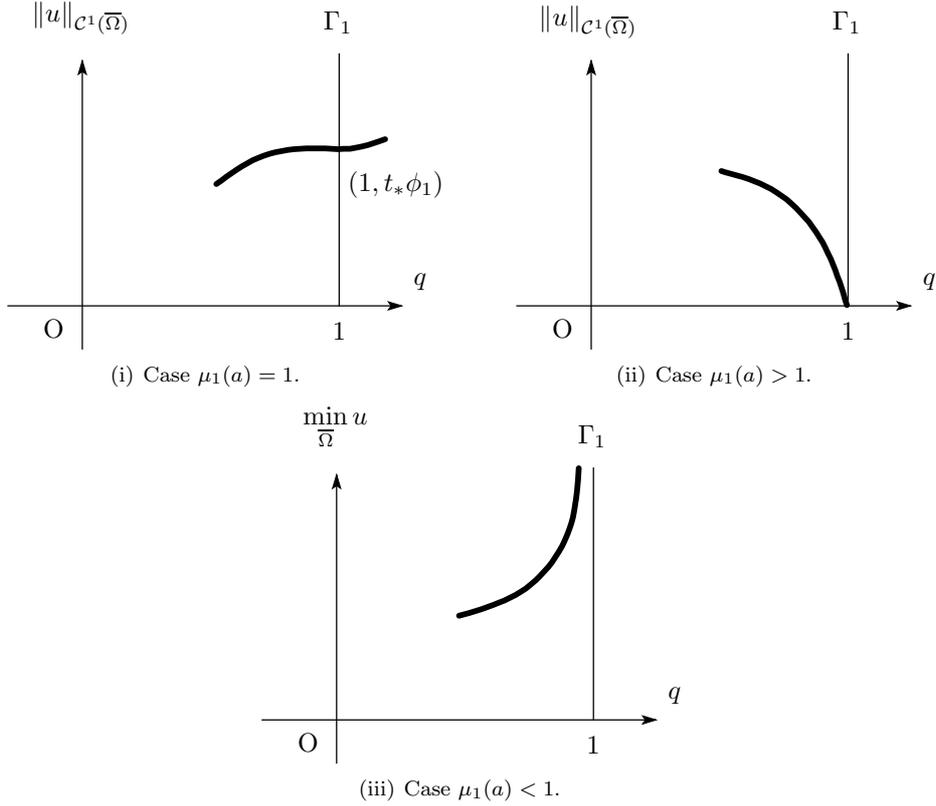
\begin{figure}[!htb] 
      \begin{center} 
    % (a)
      \subfigure[Case $\mu_1(a)=1$.]{
        %WinTpicVersion4.28b
{\unitlength 0.1in
\begin{picture}( 20.4600, 17.9500)( 19.1000,-23.8000)
% STR 2 0 3 0 Black White
% 4 4046 1949 4046 2028 5 0 0 0
% $q$
\put(40.4600,-20.2800){\makebox(0,0){$q$}}%
% VECTOR 2 0 3 0 Black White
% 2 1910 2153 3948 2153
% 
\special{pn 8}%
\special{pa 1910 2154}%
\special{pa 3948 2154}%
\special{fp}%
\special{sh 1}%
\special{pa 3948 2154}%
\special{pa 3882 2134}%
\special{pa 3896 2154}%
\special{pa 3882 2174}%
\special{pa 3948 2154}%
\special{fp}%
% LINE 2 0 3 0 Black White
% 2 3626 2153 3626 837
% 
\special{pn 8}%
\special{pa 3626 2154}%
\special{pa 3626 838}%
\special{fp}%
% STR 2 0 3 0 Black White
% 4 3626 2208 3626 2286 5 0 0 0
% $1$
\put(36.2600,-22.8600){\makebox(0,0){$1$}}%
% VECTOR 2 0 3 0 Black White
% 2 2295 2380 2295 876
% 
\special{pn 8}%
\special{pa 2296 2380}%
\special{pa 2296 876}%
\special{fp}%
\special{sh 1}%
\special{pa 2296 876}%
\special{pa 2276 944}%
\special{pa 2296 930}%
\special{pa 2316 944}%
\special{pa 2296 876}%
\special{fp}%
% STR 2 0 3 0 Black White
% 4 2148 2192 2148 2270 5 0 0 0
% O
\put(21.4800,-22.7000){\makebox(0,0){O}}%
% SPLINE 0 0 3 0 Black White
% 5 2989 1518 3213 1378 3514 1331 3696 1331 3864 1284
% 
\special{pn 30}%
\special{pa 2990 1518}%
\special{pa 3016 1500}%
\special{pa 3042 1480}%
\special{pa 3094 1444}%
\special{pa 3122 1426}%
\special{pa 3178 1394}%
\special{pa 3206 1382}%
\special{pa 3236 1370}%
\special{pa 3266 1360}%
\special{pa 3296 1352}%
\special{pa 3360 1340}%
\special{pa 3392 1336}%
\special{pa 3456 1332}%
\special{pa 3554 1332}%
\special{pa 3618 1336}%
\special{pa 3650 1334}%
\special{pa 3682 1334}%
\special{pa 3714 1328}%
\special{pa 3744 1322}%
\special{pa 3776 1314}%
\special{pa 3836 1294}%
\special{pa 3864 1284}%
\special{fp}%
% STR 2 0 3 0 Black White
% 4 3619 594 3619 673 5 0 0 0
% $\Gamma_{1}$
\put(36.1900,-6.7300){\makebox(0,0){$\Gamma_{1}$}}%
% STR 2 0 3 0 Black White
% 4 2290 550 2290 650 5 0 0 0
% $\normu$
\put(22.9000,-6.5000){\makebox(0,0){$\normu$}}%
% STR 2 0 3 0 Black White
% 4 3920 1420 3920 1520 5 0 0 0
% $\bpt$
\put(39.2000,-15.2000){\makebox(0,0){$\bpt$}}%
\end{picture}}%
          \label{fig01}
      }
      \hfill
      % (b)
      \subfigure[Case $\mu_1(a)>1$.]{
          %WinTpicVersion4.28b
{\unitlength 0.1in
\begin{picture}( 20.4600, 17.7500)( 19.1000,-23.8000)
% STR 2 0 3 0 Black White
% 4 4046 1949 4046 2028 5 0 0 0
% $q$
\put(40.4600,-20.2800){\makebox(0,0){$q$}}%
% VECTOR 2 0 3 0 Black White
% 2 1910 2153 3948 2153
% 
\special{pn 8}%
\special{pa 1910 2154}%
\special{pa 3948 2154}%
\special{fp}%
\special{sh 1}%
\special{pa 3948 2154}%
\special{pa 3882 2134}%
\special{pa 3896 2154}%
\special{pa 3882 2174}%
\special{pa 3948 2154}%
\special{fp}%
% LINE 2 0 3 0 Black White
% 2 3626 2153 3626 837
% 
\special{pn 8}%
\special{pa 3626 2154}%
\special{pa 3626 838}%
\special{fp}%
% STR 2 0 3 0 Black White
% 4 3626 2208 3626 2286 5 0 0 0
% $1$
\put(36.2600,-22.8600){\makebox(0,0){$1$}}%
% VECTOR 2 0 3 0 Black White
% 2 2295 2380 2295 876
% 
\special{pn 8}%
\special{pa 2296 2380}%
\special{pa 2296 876}%
\special{fp}%
\special{sh 1}%
\special{pa 2296 876}%
\special{pa 2276 944}%
\special{pa 2296 930}%
\special{pa 2316 944}%
\special{pa 2296 876}%
\special{fp}%
% STR 2 0 3 0 Black White
% 4 2148 2192 2148 2270 5 0 0 0
% O
\put(21.4800,-22.7000){\makebox(0,0){O}}%
% STR 2 0 3 0 Black White
% 4 3619 594 3619 673 5 0 0 0
% $\Gamma_{1}$
\put(36.1900,-6.7300){\makebox(0,0){$\Gamma_{1}$}}%
% SPLINE 0 0 3 0 Black White
% 4 3620 2150 3480 1800 3220 1540 2970 1450
% 
\special{pn 30}%
\special{pa 3620 2150}%
\special{pa 3610 2120}%
\special{pa 3602 2090}%
\special{pa 3592 2058}%
\special{pa 3582 2028}%
\special{pa 3570 1998}%
\special{pa 3560 1968}%
\special{pa 3536 1908}%
\special{pa 3522 1880}%
\special{pa 3508 1850}%
\special{pa 3494 1822}%
\special{pa 3476 1794}%
\special{pa 3460 1768}%
\special{pa 3440 1740}%
\special{pa 3422 1714}%
\special{pa 3356 1642}%
\special{pa 3308 1598}%
\special{pa 3256 1562}%
\special{pa 3228 1544}%
\special{pa 3172 1516}%
\special{pa 3112 1492}%
\special{pa 3082 1482}%
\special{pa 2988 1456}%
\special{pa 2970 1450}%
\special{fp}%
% STR 2 0 3 0 Black White
% 4 2280 570 2280 670 5 0 0 0
% $\normu$
\put(22.8000,-6.7000){\makebox(0,0){$\normu$}}%
\end{picture}}%
        \label{fig03}
      }
      \hfill
      % (c)
      \subfigure[Case $\mu_1(a)<1$.]{
          %WinTpicVersion4.28b
{\unitlength 0.1in
\begin{picture}( 20.4600, 18.1500)( 19.1000,-23.8000)
% STR 2 0 3 0 Black White
% 4 4046 1949 4046 2028 5 0 0 0
% $q$
\put(40.4600,-20.2800){\makebox(0,0){$q$}}%
% VECTOR 2 0 3 0 Black White
% 2 1910 2153 3948 2153
% 
\special{pn 8}%
\special{pa 1910 2154}%
\special{pa 3948 2154}%
\special{fp}%
\special{sh 1}%
\special{pa 3948 2154}%
\special{pa 3882 2134}%
\special{pa 3896 2154}%
\special{pa 3882 2174}%
\special{pa 3948 2154}%
\special{fp}%
% LINE 2 0 3 0 Black White
% 2 3626 2153 3626 837
% 
\special{pn 8}%
\special{pa 3626 2154}%
\special{pa 3626 838}%
\special{fp}%
% STR 2 0 3 0 Black White
% 4 3626 2208 3626 2286 5 0 0 0
% $1$
\put(36.2600,-22.8600){\makebox(0,0){$1$}}%
% VECTOR 2 0 3 0 Black White
% 2 2295 2380 2295 876
% 
\special{pn 8}%
\special{pa 2296 2380}%
\special{pa 2296 876}%
\special{fp}%
\special{sh 1}%
\special{pa 2296 876}%
\special{pa 2276 944}%
\special{pa 2296 930}%
\special{pa 2316 944}%
\special{pa 2296 876}%
\special{fp}%
% STR 2 0 3 0 Black White
% 4 2148 2192 2148 2270 5 0 0 0
% O
\put(21.4800,-22.7000){\makebox(0,0){O}}%
% STR 2 0 3 0 Black White
% 4 3619 594 3619 673 5 0 0 0
% $\Gamma_{1}$
\put(36.1900,-6.7300){\makebox(0,0){$\Gamma_{1}$}}%
% SPLINE 0 0 3 0 Black White
% 4 3550 840 3500 1160 3270 1470 2930 1610
% 
\special{pn 30}%
\special{pa 3550 840}%
\special{pa 3546 904}%
\special{pa 3540 968}%
\special{pa 3536 1000}%
\special{pa 3526 1064}%
\special{pa 3520 1096}%
\special{pa 3512 1126}%
\special{pa 3502 1156}%
\special{pa 3490 1188}%
\special{pa 3478 1216}%
\special{pa 3464 1246}%
\special{pa 3448 1274}%
\special{pa 3412 1330}%
\special{pa 3392 1356}%
\special{pa 3348 1404}%
\special{pa 3326 1426}%
\special{pa 3300 1448}%
\special{pa 3276 1468}%
\special{pa 3248 1486}%
\special{pa 3222 1502}%
\special{pa 3166 1530}%
\special{pa 3106 1554}%
\special{pa 3046 1576}%
\special{pa 2982 1596}%
\special{pa 2952 1604}%
\special{pa 2930 1610}%
\special{fp}%
% STR 2 0 3 0 Black White
% 4 2290 530 2290 630 5 0 0 0
% $\minu$
\put(22.9000,-6.3000){\makebox(0,0){$\minu$}}%
\end{picture}}%
        \label{fig02} 
      }

      \end{center}
      \caption{Bifurcating solutions in $\mathcal{P}^{\circ}$ 
	from the trivial line $\Gamma_1$.}
      \label{figbif} 
    \end{figure} 
%%%%%%%%%%%%%%%%%%%%%%%%%%%%%%%%%%%%

When $\Omega$ is a ball and $a$ is radial, we shall exhibit some
\textit{explicit} conditions on $q$ and $a$ so that $(P_{a,q})$ admits a
positive solution. This will be done via the well known sub-supersolutions
method. In Theorem \ref{cc} below we give a condition that guarantees the
existence of a positive solution (not necessarily in $\mathcal{P}^{\circ}$),
while Theorem \ref{rad2} provides us with a solution in $\mathcal{P}^{\circ}$.

Given $0<R_{0}<R$, we write%
\begin{align*}
B_{R_{0}}  &  :=\left\{  x\in\mathbb{R}^{N}:\left\vert x\right\vert
<R_{0}\right\}  ,\\
A_{R_{0},R}  &  :=\left\{  x\in\mathbb{R}^{N}:R_{0}<\left\vert x\right\vert
<R\right\}  ,\\
\omega_{N-1}  &  :=\text{surface area of the unit sphere }\partial B_{1}\text{
in }\mathbb{R}^{N}.
\end{align*}

If $f$ is a radial function, we write (with a slight abuse of notation)
$f\left(  x\right)  :=f\left(  \left\vert x\right\vert \right)  :=f\left(
r\right)  $. We first consider the case that $\text{supp }a^{+}$ is contained
in $B_{R_{0}}$ for some $R_{0}\in(0,R)$.
%\marginpar{modified}

\begin{theorem}
\label{cc}
%\marginpar{\textit{modified}}
Let $\Omega:=B_{R}$ and $a\in L^{\infty}\left(  \Omega\right)  $\textit{ be a
radial function such that }$\int_{\Omega}a<0$.
%\marginpar{\textit{modified}}
Assume that there exists $R_{0}>0$ such that:
%\marginpar{\it modified}

\begin{itemize}
\item $a \geq0$ in $B_{R_{0}}$;

\item $a \leq0$ in $A_{R_{0},R}$;

\item $r\rightarrow a(r)$ is differentiable and nonincreasing in $\left(
R_{0},R\right)  $, and%
\begin{equation}
\frac{1-q}{1+q}\int_{A_{R_{0},R}}a^{-}\leq\int_{B_{R_{0}}}a^{+}%
.\label{inferno}%
\end{equation}
Then $(P_{a,q})$ has a positive
%\marginpar{\textit{modified}}
solution, 
	% \marginpar{\textit{modified}} 
which is unique if $(H_{1})$ holds.
\end{itemize}
\end{theorem}

\begin{remark}\strut 
		% \marginpar{\it modified}
\begin{enumerate}
\item The condition \eqref{inferno} can also be formulated as
\begin{equation}\label{cq}
\frac{-\int_{\Omega}a}{\int_{\Omega}\left\vert a\right\vert }\leq q<1.
\end{equation}
In particular, we see that \eqref{inferno} is satisfied if $q$ is close
enough to $1$. Note that if we 
replace $a$ by 
\[
a_\delta =a^+ - \delta a^-,  \ \mbox{ with } \ 
\delta > \delta_0 := \frac{\int_\Omega a^+}{\int_\Omega a^-}, 
\]
% keep $a^+$ fixed and make $a^-$ large 
then the left-hand side in \eqref{cq} approaches $1$ as $\delta \to \infty$, 
so that this condition becomes very restrictive for $a_\delta$ as $\delta \to 
\infty$. On the other side,  % \marginpar{\textit{modified}}, 
we have that $\int_\Omega a_\delta \to 0^-$ as $\delta \to \delta_0^+$, so that
\eqref{cq} becomes much less constraining for $a_\delta$ as $\delta \to \delta_0^+$. 
A similar argument will be used in Remark \ref{rem:I01}. 

\item As one can see from the proof of Theorem \ref{cc}, the condition \eqref{inferno} guarantees the existence of a positive subsolution for the
corresponding Dirichlet problem. Thus, since arbitrarily large supersolutions
can be easily obtained in the Dirichlet case (see e.g. \cite[Remark
1.1]{nodea}), it follows that (\ref{inferno}) ensures the existence of a
positive solution for the analogous Dirichlet problem. Moreover, we point out
that this condition substantially improves some of the results known in that
case (see \cite[Section 3]{nodea}).
\end{enumerate}
\end{remark}

Next we consider the case that $\text{supp }a^{-}$ is contained in $B_{R_{0}}$
for some $R_{0}\in(0,R)$.

\begin{theorem}
\label{rad2}
%\marginpar{\textit{modified}}
Let $\Omega:=B_{R}$ and $a\in L^{\infty}\left(  \Omega\right)  $\textit{ be a
radial function such that }$\int_{\Omega}a<0$. Assume that there exists
$R_{0}\in\left(  0,R\right)  $ such that $a\geq0$ in $A_{R_{0},R}$, and
%\marginpar{\textit{modified (1.3)}}%
\begin{equation}
\frac{1-q}{2q+N\left(  1-q\right)  }\omega_{N-1}R_{0}^{N}\left\Vert
a^{-}\right\Vert _{L^{\infty}(B_{R_{0}})}<\int_{A_{R_{0},R}}a^{+}.
\label{sipi}%
\end{equation}
Then $(P_{a,q})$ has a unique solution $u\in\mathcal{P}^{\circ}$.
\end{theorem}

\begin{remark} 
Observe that unlike Theorem \ref{cc}, no differentiability nor 
monoto-nicity condition is imposed on $a^{-}$ in Theorem \ref{rad2}. 
Note again that
(\ref{sipi}) is also clearly satisfied if $q$ is close enough to $1$.
\end{remark}

Our next results concern the sets
\[
\mathcal{A}=\mathcal{A}_{a}:=\{q\in(0,1):\text{any nontrivial nonnegative
solution of }(P_{a,q})\text{ lies in }\mathcal{P}^{\circ}\}
\]
and
\[
\mathcal{I}=\mathcal{I}_{a}:=\{q\in(0,1):(P_{a,q})\text{ has a solution }%
u\in\mathcal{P}^{\circ}\}.
\]

	%\marginpar{\textit{modified}}
We observe that if $(H_{0})$ holds then
$(P_{a,q})$ has a nontrivial nonnegative solution for any $0<q<1$ (see e.g.
the proof of \cite[Corollary 1.8]{krqu}), so that $\mathcal{A}\subseteq
\mathcal{I}$. In \cite[Theorem 1.9]{krqu}, we proved, under $(H_{0})$ and
$(H_{1})$, that $\mathcal{A}=(q_{a},1)$ for some $q_{a}\in\lbrack0,1)$.

Let us now introduce the following assumptions:
%\marginpar{\textit{modified}}%
%\marginpar{\textit{modified}}%
\[
\Omega_{+}\text{ is connected and }\left\vert (\text{supp }a^{+}%
)\setminus\Omega_{+}\right\vert =0,\leqno{(H_1')}
\]%
\[
\partial\Omega_{+}\text{ satisfies the inner sphere condition with respect to
}\Omega_{+}.\leqno{(H_+)}
\]
Note
	% \marginpar{\textit{modified}} 
that $(H_{1}^{\prime})$ corresponds to
$(H_{1})$ with $\Omega_{+}$ consisting of a single connected component.

\begin{remark}
\label{h}If $a$ is H\"{o}lder continuous, $\Omega_{+}$ has finitely
%\marginpar{\textit{modified}}
many connected components and $(H_{+})$ holds, then \cite[Theorem 3.1]{BPT} shows,
in particular, that $(P_{a,q})$ has at most one nonnegative solution $u$ such
that $u>0$ in $\Omega_{+}$. Let us observe that their proof is still valid for
$a\in L^{r}(\Omega)$ with $r>N$, assuming now $(H_{1})$ and $(H_{+})$.
\end{remark}

We shall complement \cite[Theorem 1.9]{krqu} as follows:

\begin{theorem}
\strut\label{ti}

\begin{enumerate}
\item[(i)] Assume $(H_{0})$. Then $\mathcal{I}\not =\emptyset$. Moreover:

\begin{enumerate}
\item[(i1)] Either $\mathcal{I}=[q_{i},1)$ with $q_{i}\in\left(  0,1\right)  $
or $\mathcal{I}=\left(  q_{i},1\right)  $ with $q_{i}\in\lbrack0,1)$.

\item[(i2)] If $(H_{1}^{\prime})$ and $(H_{+})$ hold, then for all
$q\in\left(  0,1\right)  $,
%\marginpar{\textit{modified}}
there exists a unique nontrivial nonnegative solution of $(P_{a,q})$. In
particular, $\mathcal{I}=\mathcal{A}$.
%\marginpar{modified}

\end{enumerate}

\item[(ii)] Let $\Omega:=(x_{0},x_{1})\subset\mathbb{R}$. Given $q\in\left(
0,1\right)  $, there exists $a\in\mathcal{C}(\overline{\Omega})$ such that
$q\in\mathcal{I}\setminus\mathcal{A}$.
\end{enumerate}
\end{theorem}

%Under $(H_{1})$, it is already known from our previous work \cite[Corollary 1.8]{krqu} that there exists $q_{0}\in(0,1)$ such that $(P_{a,q})$ has a solution in $\mathcal{P}^{\circ}$ (which is the unique nontrivial nonnegative solution of $(P_{a,q})$) for $q\in(q_{0},1)$. In this sense, Corollary \ref{c1} improves \cite[Corollary 1.8]{krqu}. Furthermore, the above corollary provides the asymptotic behavior of $u_{q}$ as $q\rightarrow1^{-}$.\newline

\begin{remark}
\label{q0} \strut
\begin{enumerate}
\item Let $q_{1}\in\left(  0,1\right)  $, and define $\Omega:=\left(  0,\pi\right)
$,
\[
r:=\frac{2}{1-q_{1}}\in\left(  2,\infty\right),\quad \text{and} \quad a\left(  x\right)
:=r^{1-\frac{2}{r}}\left(  1-r\cos^{2}x\right)  ,\quad\text{for }x\in\Omega.
\]
One can check that $u\left(  x\right)  :=\frac{\sin^{r}x}{r}\ $is a (strictly
positive in $\Omega$) solution of 
\[
\left\{
\begin{array}
[c]{lll}%
-u^{\prime\prime}=a(x)u^{q_{1}} & \mathrm{in} & \Omega,\\
u^{\prime}=0 & \mathrm{on} & \partial\Omega.
\end{array}
\right.
\]
	% \marginpar{\textit{modified}} 
It follows that $q_{1}\not \in \mathcal{A}$
because  $u\not \in \mathcal{P}^{\circ}$. Now, since $a$ satisfies
$(H_{1}^{\prime})$ and $(H_{+})$, we deduce from Theorem \ref{ti} (i2) that 
$u$ is the unique nontrivial nonnegative solution of $(P_{a,q_{1}})$, and 
$\mathcal{I}=\mathcal{A}$. Consequently, we have 
$q_{1}\not \in \mathcal{I}$. In particular, if $q_{0}\left(  a\right)  $ 
is given by Theorem \ref{c1}, then
$q_{0}\left(  a\right)  \geq q_{1}$. In the same way, if $q_{i}$ is provided
by Theorem 
\ref{ti} (i1), then $q_{i}\geq q_{1}$.\\

\item After the corresponding modifications, Theorem \ref{ti} (i) holds also for the Dirichlet counterpart of $(P_{a,q})$, replacing $(H_0)$ by the condition that $a^+ \not \equiv 0$. As a matter of fact, in this case one can check that the proof of Theorem \ref{ti} (i) can be carried out, using now \cite[Theorems 2.1, 2.2, and Lemma 2.3]{bandle}. 
\end{enumerate}
\end{remark}

Finally, we shall % \marginpar{\textit{modified}} 
investigate the existence of nonnegative \textit{dead core}
solutions of $(P_{a,q})$. Following \cite{bandle, BPT}, the set $\{x\in
\Omega:u(x)=0\}$ is called the \textit{dead core} of a nontrivial nonnegative
solution $u$ of $(P_{a,q})$. Let us mention that in the proof of Theorem
\ref{ti} (ii) we shall see that, when $N=1$, for any $q\in(0,1)$ there exists
$a$ with $(P_{a,q})$ admitting a solution in $\mathcal{P}^{\circ}$ and also
nonnegative solutions with nonempty dead cores.
%\marginpar{\textit{modified}}

	% \marginpar{\textit{modified}}
Next we give some \textit{sufficient} conditions
for the existence of dead core solutions of $(P_{a,q})$. We introduce the
following condition:
%\marginpar{\textit{modified}}
%\marginpar{\textit{modified }$H_{2}$ \textit{in order to be able to
%use Rem. 3.1 in the proof of Th. 1.10 (ii)}}%
\[
b_{1}\in L^{\infty}(\Omega),b_{2}\in\mathcal{C}(\overline{\Omega}),b_{1}%
,b_{2}\geq0\text{ and }\text{supp }b_{1}\cap\{x\in\Omega:b_{2}%
(x)>0\}=\emptyset.\leqno{(H_2)}
\]
%\marginpar{\textit{modified}}
Given a nonempty open subset $G\subseteq\Omega$ and $\sigma>0$, we set
\begin{equation}
G_{\sigma}:=\{x\in G:\text{dist}(x,\partial G)>\sigma\}.\label{sig}%
\end{equation}
We call the set $\Omega\setminus\overline{\Omega_{\sigma}}$ a \textit{tubular
neighborhood} of $\partial\Omega$.

\begin{theorem}
\label{dc} \strut

\begin{enumerate}
\item Let $q\in(0,1)$, and assume that $(H_{1}^{\prime})$ holds and
$\Omega_{+}$ contains a tubular neighborhood of $\partial\Omega$. Then,
\textit{every} nontrivial nonnegative solution of $(P_{a,q})$ is positive on
$\partial\Omega$. In particular, if $u$ is a nontrivial nonnegative solution
of $(P_{a,q})$, then either $u\in\mathcal{P}^{\circ}$ or $u$ has a nonempty
dead core.

\item Let $a_{\delta}:=b_{1}-\delta b_{2}$, with $b_{1},b_{2}\not \equiv 0$
%\marginpar{\textit{modified}}
satisfying $(H_{2})$, and $\delta> 0$. If we set $G:=\{x\in\Omega
:b_{2}(x)>0\}$ then, given $0<\overline{q}<1$ and $\sigma>0$, there exists
$\delta_{0}=\delta_{0}(\sigma,\overline{q})>0$ such that any nontrivial
nonnegative solution of $(P_{a_{\delta},q})$ with $q\in(0,\overline{q}]$
vanishes in $G_{\sigma}$ if $\delta\geq\delta_{0}$.
\end{enumerate}
\end{theorem}

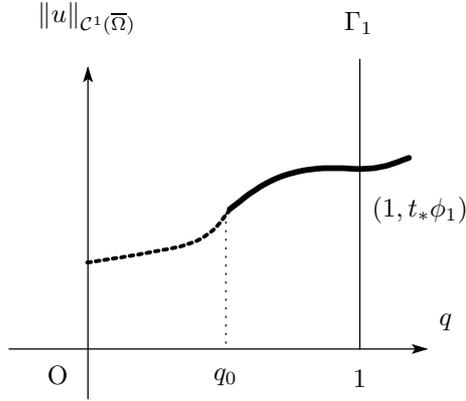
\begin{figure}[!htb]
  	   \begin{center}
		%WinTpicVersion4.28b
{\unitlength 0.1in
\begin{picture}( 21.7100, 20.5200)( 19.1000,-26.5000)
% STR 2 0 3 0 Black White
% 4 4171 2155 4171 2246 5 0 0 0
% $q$
\put(41.7100,-22.4600){\makebox(0,0){$q$}}%
% VECTOR 2 0 3 0 Black White
% 2 1910 2389 4068 2389
% 
\special{pn 8}%
\special{pa 1910 2390}%
\special{pa 4068 2390}%
\special{fp}%
\special{sh 1}%
\special{pa 4068 2390}%
\special{pa 4002 2370}%
\special{pa 4016 2390}%
\special{pa 4002 2410}%
\special{pa 4068 2390}%
\special{fp}%
% LINE 2 0 3 0 Black White
% 2 3726 2389 3726 878
% 
\special{pn 8}%
\special{pa 3726 2390}%
\special{pa 3726 878}%
\special{fp}%
% STR 2 0 3 0 Black White
% 4 3726 2453 3726 2542 5 0 0 0
% $1$
\put(37.2600,-25.4200){\makebox(0,0){$1$}}%
% VECTOR 2 0 3 0 Black White
% 2 2318 2650 2318 923
% 
\special{pn 8}%
\special{pa 2318 2650}%
\special{pa 2318 924}%
\special{fp}%
\special{sh 1}%
\special{pa 2318 924}%
\special{pa 2298 990}%
\special{pa 2318 976}%
\special{pa 2338 990}%
\special{pa 2318 924}%
\special{fp}%
% STR 2 0 3 0 Black White
% 4 2162 2434 2162 2524 5 0 0 0
% O
\put(21.6200,-25.2400){\makebox(0,0){O}}%
% SPLINE 0 0 3 0 Black White
% 5 3052 1660 3290 1500 3608 1445 3801 1445 3979 1391
% 
\special{pn 30}%
\special{pa 3052 1660}%
\special{pa 3104 1620}%
\special{pa 3128 1600}%
\special{pa 3154 1580}%
\special{pa 3182 1562}%
\special{pa 3208 1544}%
\special{pa 3236 1528}%
\special{pa 3262 1514}%
\special{pa 3292 1500}%
\special{pa 3322 1488}%
\special{pa 3352 1478}%
\special{pa 3382 1470}%
\special{pa 3414 1462}%
\special{pa 3444 1456}%
\special{pa 3476 1452}%
\special{pa 3510 1448}%
\special{pa 3542 1446}%
\special{pa 3640 1446}%
\special{pa 3704 1450}%
\special{pa 3736 1450}%
\special{pa 3800 1446}%
\special{pa 3830 1440}%
\special{pa 3862 1434}%
\special{pa 3922 1414}%
\special{pa 3952 1402}%
\special{pa 3980 1392}%
\special{fp}%
% STR 2 0 3 0 Black White
% 4 3719 599 3719 689 5 0 0 0
% $\Gamma_{1}$
\put(37.1900,-6.8900){\makebox(0,0){$\Gamma_{1}$}}%
% STR 2 0 3 0 Black White
% 4 2313 548 2313 663 5 0 0 0
% $\normu$
\put(23.1300,-6.6300){\makebox(0,0){$\normu$}}%
% STR 2 0 3 0 Black White
% 4 4038 1547 4038 1662 5 0 0 0
% $\bpt$
\put(40.3800,-16.6200){\makebox(0,0){$\bpt$}}%
% SPLINE 0 2 3 0 Black White
% 4 3053 1662 2926 1800 2725 1869 2313 1938
% 
\special{pn 30}%
\special{pn 20}%
\special{pa 3054 1662}%
\special{pa 3042 1678}%
\special{fp}%
\special{pa 3027 1697}%
\special{pa 3015 1713}%
\special{fp}%
\special{pa 3000 1733}%
\special{pa 2987 1748}%
\special{fp}%
\special{pa 2970 1766}%
\special{pa 2955 1778}%
\special{fp}%
\special{pa 2935 1793}%
\special{pa 2919 1804}%
\special{fp}%
\special{pa 2898 1817}%
\special{pa 2880 1826}%
\special{fp}%
\special{pa 2857 1835}%
\special{pa 2838 1841}%
\special{fp}%
\special{pa 2815 1849}%
\special{pa 2796 1855}%
\special{fp}%
\special{pa 2772 1860}%
\special{pa 2753 1864}%
\special{fp}%
\special{pa 2728 1869}%
\special{pa 2709 1872}%
\special{fp}%
\special{pa 2684 1877}%
\special{pa 2665 1880}%
\special{fp}%
\special{pa 2641 1885}%
\special{pa 2621 1889}%
\special{fp}%
\special{pa 2597 1893}%
\special{pa 2578 1897}%
\special{fp}%
\special{pa 2553 1901}%
\special{pa 2534 1904}%
\special{fp}%
\special{pa 2509 1908}%
\special{pa 2490 1912}%
\special{fp}%
\special{pa 2465 1915}%
\special{pa 2446 1918}%
\special{fp}%
\special{pa 2422 1922}%
\special{pa 2402 1925}%
\special{fp}%
\special{pa 2378 1929}%
\special{pa 2358 1932}%
\special{fp}%
\special{pa 2334 1936}%
\special{pa 2314 1938}%
\special{fp}%
% LINE 2 2 3 0 Black White
% 2 3031 1682 3031 2382
% 
\special{pn 8}%
\special{pa 3032 1682}%
\special{pa 3032 2382}%
\special{dt 0.045}%
% STR 2 0 3 0 Black White
% 4 3031 2444 3031 2535 5 0 0 0
% $q_0$
\put(30.3100,-25.3500){\makebox(0,0){$q_0$}}%
\end{picture}}%
	  \caption{A possible bifurcation curve of the unique 
	nontrivial nonnegative solution when $\mu_1(a)=1$, conditions $(H_1'), (H_+)$ hold, and $\Omega_+$ contains a tubular neighborhood of $\partial \Omega$. Here the full curve represents solutions in $\mathcal{P}^{\circ}$, whereas the dotted curve represents dead core solutions.} 
	\label{fig04} 
	  \end{center}
	    \end{figure}

\begin{remark}
\label{rem:dc} \strut

\begin{enumerate}
\item The conclusion of Theorem \ref{dc} (ii) still holds if $a_{\delta
}:=b-\delta\chi_{G}$, with $b,G$ satisfying
\[
\left\{
\begin{array}
[c]{l}%
b\in L^{\infty}(\Omega),\ 0\not \equiv b\geq0,\mbox{ and }\\
\emptyset\not =G\subset\Omega\mbox{ is an open subset such that }\mathrm{supp}%
\ b\cap G=\emptyset.
\end{array}
\right.  \leqno{(H_2')}
\]
Here $\delta>0$ and $\chi_{G}$ is the characteristic function of $G$.

\item Let $a_{\delta}:= b_{1} - \delta b_{2}$ with $b_{1},b_{2} \not \equiv 0$
satisfying $(H_{2})$, and $\delta> 0$.

\begin{enumerate}
\item[(ii1)] In addition to $(H_{1}^{\prime})$, let us 
	% \marginpar{\textit{modified}} 
assume that
\[
\text{supp }b_{1}\cup\{x\in\Omega:b_{2}(x)>0\}=\Omega. 
\]
Let 
$q \in (0,1)$. Theorem \ref{dc} (ii) then shows that the support of any
nontrivial nonnegative solution of $(P_{a_{\delta}, q})$ approaches $\Omega_{+}$ (in some
sense) as $\delta\rightarrow\infty$.

\item[(ii2)] Combining Theorem \ref{c1} and
%Proposition \ref{prop:dead} or
Theorem \ref{dc} (ii), we find $\delta_{1}>0$ and $0<q_{1}\leq q_{0}<1$ such
that any nontrivial nonnegative solution of $(P_{a_{\delta},q})$ with
$\delta=\delta_{1}$ has a nonempty dead core for $q\in(0,q_{1}]$, whereas this
problem has a unique solution in $\mathcal{P}^{\circ}$ and no other nontrivial
nonnegative solutions for $q\in(q_{0},1)$. Furthermore, according to Theorem
\ref{ti} (i2) and Theorem \ref{dc} (i), we see that if $(H_{1}^{\prime})$ and
$(H_{+})$ hold
%\marginpar{\textit{modified}}
and $\Omega_{+}$ contains a tubular neighborhood of $\partial\Omega$, then
$q_{1}=q_{0}$, and the nontrivial nonnegative solution for $q\in(0,q_{0}]$ is
also unique 
(see Figure \ref{fig04}).

\item[(ii3)] As we shall see from its proof, Theorem \ref{dc} (ii) holds also 
for the Dirichlet counterpart of $(P_{a_{\delta},q})$. In particular, it
complements \cite[Theorem 1.1]{krqu} as follows: given $q\in(0,1)$ there exist
$0<\delta_{1}<\delta_{0}$ such that every nontrivial nonnegative solution 
$u$ of 
\[
\left\{ 
\begin{array}
[c]{lll}%
-\Delta u=a_{\delta}(x)u^{q} & \mathrm{in} & \Omega,\\
u=0 & \mathrm{on} & \partial\Omega,
\end{array}
\right.
\]
satisfies $u>0$ in $\Omega$ and $\frac{\partial u}{\partial\nu}<0$ on
$\partial\Omega$ for $\delta<\delta_{1}$, whereas $u$ has a nonempty dead core
for $\delta>\delta_{0}$. \newline
\end{enumerate}
\end{enumerate}
\end{remark}

%%%%%%%%%%%%%%%%%%%%%%%%%%%%%%%%%%%%%%%%
%	% \begin{figure}[H] 
%%%%%%%%%%%%%%%%%%%%%%%%%%%%%%%%%%%%%%%%%%%

The rest of the paper is organized as follows: in Section $2$ we establish
some bifurcation results and stability properties for solutions in
$\mathcal{P}^{\circ}$ of $(P_{a,q})$, whereas Section $3$ is devoted to the
proof of Theorems \ref{c1}, \ref{cc} and \ref{rad2}. In Section $4$ we prove
Theorem \ref{ti} and some corollaries of it.
%\marginpar{\textit{modified}}
Finally, Section \ref{sec:core} is concerned with the existence of dead core
solutions and the proof of Theorem \ref{dc}.
%\marginpar{\textit{modified}}
\newline

%%%

\section{Bifurcation analysis}

This section is devoted to the bifurcation analysis of $(P_{a,q})$, where $q$
is the bifurcation parameter. First we establish, under $(H_{0})$ and
$(H_{1})$, some \textit{a priori} bounds for nontrivial nonnegative solutions
of $(P_{a,q})$, which imply that no nontrivial nonnegative solutions bifurcate
from zero or from infinity at any $q\in\lbrack0,1)$. More precisely, we shall
see that given $q\in\lbrack0,1)$ there exists no sequence $\left\{
q_{n}\right\}  \subset(0,1)$ such that $q_{n} \to q$ and $(P_{a,q_{n}})$ has a
nontrivial nonnegative solution $u_{n}$ satisfying $u_{n}\rightarrow0$ in
$\mathcal{C}(\overline{\Omega})$ or $\Vert u_{n}\Vert_{\infty}\rightarrow
\infty$.

\begin{proposition}
\label{t2.1} \strut

\begin{enumerate}
\item Assume $(H_{1})$. Then, given $q_{1}\in(0,1)$, there exists $C>1$ such
that $\Vert u\Vert_{L^{\infty}(\Omega)}>C^{-1}$ for all nontrivial nonnegative
solutions of $(P_{a,q})$ with $q\in(0,q_{1}]$.

\item Assume $(H_{0})$. Then, given $q_{1}\in(0,1)$, there exists $C>1$ such
that $\Vert u\Vert_{L^{\infty}(\Omega)}<C$ for all nontrivial nonnegative
solutions of $(P_{a,q})$ with $q\in(0,q_{1}]$. \newline
\end{enumerate}
\end{proposition}

\textbf{Proof.}

\begin{enumerate}
\item First we obtain an \textit{a priori} bound from below. Assume by
contradiction that there exist $0<q_{n}\leq\overline{q}<1$ and $u_{n}$
nontrivial nonnegative solutions of $(P_{a,q_{n}})$ such that $u_{n}%
\rightarrow0$ in $\mathcal{C}(\overline{\Omega})$. Then, thanks to $(H_{1})$,
we may assume that $u_{n}\not \equiv 0$ in some fixed subdomain $\Omega
^{\prime}\subset\Omega_{+}$. By the strong maximum principle, we deduce that
$u_{n}>0$ in $\Omega^{\prime}$.

We fix $c>0$ sufficiently large such that $\lambda_{1}(ca,\Omega^{\prime})<1$,
where $\lambda_{1}(m,\Omega)$ denotes the first positive eigenvalue of the
Dirichlet problem
\[
\left\{
\begin{array}
[c]{lll}%
-\Delta\phi=\lambda m\left(  x\right)  \phi & \mathrm{in} & \Omega,\\
\phi=0 & \mathrm{on} & \partial\Omega,
\end{array}
\right.
\]
and observe that $v_{n}:=c^{\frac{1}{1-q_{n}}}u_{n}$ are nontrivial
nonnegative solutions of $(P_{ca,q_{n}})$. We now apply \cite[Lemma 2.5]{krqu}
to get a ball $B\subset\Omega^{\prime}$ and a positive function $\psi$ in $B$
such that
\[
v_{n}=c^{\frac{1}{1-q_{n}}}u_{n}\geq\psi\quad\mbox{ in }B,
\]
where $\psi$ and $B$ do not depend on $n$. It follows that $u_{n}\geq
c^{-\frac{1}{1-q_{n}}}\psi$ in $B$, which provides a contradiction, since
$q_{n}\leq\overline{q}<1$. \newline

\item We obtain now an \textit{a priori} bound from above. Assume to the
contrary that there exist $0<q_{n}\leq\overline{q}<1$ and $u_{n}$ nontrivial
nonnegative solutions of $(P_{a,q_{n}})$ such that $\Vert u_{n}\Vert:=\Vert
u_{n}\Vert_{H^{1}(\Omega)}\rightarrow\infty$. We set $v_{n}:=\frac{u_{n}%
}{\Vert u_{n}\Vert}$, so that we may assume that $v_{n}\rightharpoonup v_{0}$
in $H^{1}(\Omega)$ and $v_{n}\rightarrow v_{0}$ in $L^{s}(\Omega)$ for
$s\in\lbrack1,2^{\ast})$. From $(P_{a,q_{n}})$ we have that
\[
\int_{\Omega}|\nabla v_{n}|^{2}=\left(  \int_{\Omega}a(x)v_{n}^{q_{n}%
+1}\right)  \Vert u_{n}\Vert^{q_{n}-1}.
\]
Since $q_{n}\leq\overline{q}<1$, it follows that $\int_{\Omega}|\nabla
v_{n}|^{2}\rightarrow0$. Hence, we deduce that $v_{n}\rightarrow v_{0}$ in
$H^{1}(\Omega)$, and $v_{0}$ is a positive constant. Finally, since
$\int_{\Omega}a(x)v_{n}^{q_{n}+1}>0$ we derive that $\int_{\Omega}a\geq0$,
which contradicts $(H_{0})$. By elliptic regularity, we have the desired
conclusion.
%The proof is complete.
$\blacksquare$\newline
\end{enumerate}

In view of Proposition \ref{t2.1}, we see that, under $(H_{0})$ and $(H_{1})$,
	% \marginpar{\textit{modified}} 
bifurcation from zero or from infinity can only
occur at $q=1$. As already mentioned, we shall look at $q$ as a bifurcation
parameter in $(P_{a,q})$, and then
%\marginpar{\textit{modified}}
seek for bifurcating solutions in $\mathcal{P}^{\circ}$ from the trivial line
$\Gamma_{1}=\{(1,t\phi_{1}):t>0\}$ when $\mu_{1}(a)=1$. To this end, we employ
the Lyapunov-Schmidt reduction for $(P_{a,q})$, based on the positive
eigenfunction $\phi_{1}$. We set
\[
A:=-\Delta-a(x)\quad\mbox{and}\quad D(A):=\left\{  u\in W^{2,r}(\Omega
):\frac{\partial u}{\partial\nu}=0\mbox{ on }\partial\Omega\right\}  .
\]
The usual decomposition of $D(A)$ is given by the formula
\[
D(A)=\mathrm{Ker}A+X_{2};\quad u=t\phi_{1}+w,
\]
where $t=\int_{\Omega}u\phi_{1}$, and $w=u-(\int_{\Omega}u\phi_{1})\phi_{1}$.
So, $X_{2}$ is characterized as
\[
X_{2}=\left\{  w\in W^{2,r}(\Omega):\int_{\Omega}w\phi_{1}=0\right\}  .
\]
	%\marginpar{\textit{modified}} 
On the other hand, put $Y:=L^{r}(\Omega
)=Y_{1}+R(A)$, where
\[
R(A):=\left\{  f\in L^{r}(\Omega):\int_{\Omega}f\phi_{1}=0\right\},
\]
and  
	% \marginpar{\textit{modified}} 
$$ Y_1 = \langle \phi_1 \rangle := \{ s \phi_1 : s \in \mathbb{R} \}. $$  
%$$Y_1=\left\{ \left(\int_{\Omega}f\phi_{1}\right) \phi_1 : f\in L^{r}(\Omega)\right\}.$$
Let $Q$ be the projection of $Y$ to $R(A)$, given by
\[
Q[f]:=f-\left(  \int_{\Omega}f\phi_{1}\right)  \phi_{1}.
\]
We reduce $(P_{a,q})$ to the following coupled equations:
\begin{align*}
&  Q[Au]=Q[a\left(  x\right)  \left(  u^{q}-u\right)  ],\\
&  (1-Q)[Au]=(1-Q)[a\left(  x\right)  \left(  u^{q}-u\right)  ].
\end{align*}
The first equation yields
\begin{equation}
-\Delta w-a(x)w=Q[a\left(  x\right)  \{(t\phi_{1}+w)^{q}-(t\phi_{1}%
+w)\}],\label{beq01}%
\end{equation}
where we have used the fact that $\int_{\Omega}Au\phi_{1}=\int_{\Omega}%
uA\phi_{1}=0$. The second equation implies that
\begin{align*}
0 &  =(1-Q)[a\left(  x\right)  (u^{q}-u)]\\
&  =\left(  \int_{\Omega}a(x)\left\{  (t\phi_{1}+w)^{q}-(t\phi_{1}+w)\right\}
\phi_{1}\right)  \phi_{1},
\end{align*}
and thus, that
\begin{equation}
0=\int_{\Omega}a(x)\left\{  (t\phi_{1}+w)^{q}-(t\phi_{1}+w)\right\}  \phi
_{1}.\label{beq02}%
\end{equation}

Now, we see that $(q,t,w)=(1,t,0)$ satisfies (\ref{beq01}) and (\ref{beq02})
for any $t>0$. So, first we solve \eqref{beq01} with respect to $w$, around
$(q,t,w)=(1,t_{0},0)$ for a fixed $t_{0}>0$. To this end, we introduce the
mapping $F:(1-\delta,1+\delta)\times(t_{0}-d, t_{0}+d)\times B_{\rho
}(0)\rightarrow R(A)$ given by
\[
F(q,t,w):=-\Delta w-a\left(  x\right)  w-Q[a\left(  x\right)  \{(t\phi
_{1}+w)^{q}-(t\phi_{1}+w)\}],
\]
where $B_{\rho}(w)$ is the ball in $X_{2}$ centered at $w$ and with radius
$\rho> 0$. It is clear that $F(1,t_{0},0)=0$. Moreover, the Fr\'{e}chet
derivative $F_{w}(q,t,w): X_{2}\rightarrow R(A)$ is given by
\[
F_{w}(q,t,w)\varphi=-\Delta\varphi-a\left(  x\right)  \varphi-Q[a\left(
x\right)  (q(t\phi_{1}+w)^{q-1}-1)\varphi].
\]
We see that $F_{w}(1,t_{0},0)\varphi=-\Delta\varphi-a\left(  x\right)
\varphi$. Hence,
\[
F_{w}(1,t_{0},0)\varphi=0\Longleftrightarrow\varphi=c\phi_{1}%
\ \mbox{ for some }c>0.
\]
Since $\varphi\in X_{2}$, it follows that $\int_{\Omega}(c\phi_{1})\phi_{1}%
=0$, and thus $c=0$. This means that $F_{w}(1,t_{0},0)$ is injective. It is
also surjective from the fact that $\int_{\Omega}f\phi_{1}=0$ if and only if
there exists $\varphi$ such that%
\[
\left\{
\begin{array}
[c]{lll}%
-\Delta\varphi-a(x)\varphi=f & \mathrm{in} & \Omega,\\
\frac{\partial\varphi}{\partial\nu}=0 & \mathrm{on} & \partial\Omega.
\end{array}
\right.
\]
Since $F_{w}(1,t_{0},0)$ is continuous, from the Bounded Inverse Theorem we
infer that $F_{w}(1,t_{0},0)$ is an isomorphism. Hence, the implicit function
theorem applies, and consequently, we have
\begin{align*}
&  F(q,t,w)=0,\ (q,t,w)\simeq(1,t_{0},0)\\
&  \Longleftrightarrow w=w(q,t),\ (q,t)\simeq(1,t_{0}%
)\mbox{ such that }w(1,t_{0})=0.
\end{align*}
We plug $w(q,t)$ into \eqref{beq02} to get the following bifurcation equation
in $\mathbb{R}^{2}$:
\[
\Phi(q,t):=\int_{\Omega}a(x)\{(t\phi_{1}+w(q,t))^{q}-(t\phi_{1}+w(q,t)\}\phi
_{1}=0,\quad(q,t)\simeq(1,t_{0}).
\]

We are now in position to prove the following result:
%\marginpar{\textit{added}}

\begin{theorem}
\label{t1} Assume $(H_{0})$. If $\mu_{1}(a) = 1$, then the following
assertions hold:

\begin{enumerate}
\item Assume that $(q_{n},u_{n})\in(0,1)\times\mathcal{P}^{\circ}$ are
solutions of $(P_{a,q_{n}})$ such that $(q_{n},u_{n})\rightarrow(1,t\phi
_{1})\in\Gamma_{1}$ in $\mathbb{R}\times W^{2,r}(\Omega)$ for some $t>0$.
Then, we have $t=t_{*}$, where $t_{*}$ is given by \eqref{et*}.

\item The set of positive solutions of $(P_{a,q})$ consists of $\Gamma_{1}%
\cup\Gamma_{2}$ in a neighborhood of $(q,u)=(1,t_{\ast}\phi_{1})$ in
$\mathbb{R}\times W^{2,r}(\Omega)$, where
\[
\Gamma_{2}:=\{(q,\,t(q)\phi_{1}+w(q,t(q))):|q-1|<\delta_{\ast}\}
\quad\mbox{for some $\delta_{\ast} > 0$.}
\]
Here $t(q)$ and $w(q,t(q))$ are smooth with respect to $q$ and satisfy
$t(1)=t_{\ast}$ and $w(1,t_{\ast})=0$.
\end{enumerate}
\end{theorem}

%%%

\textbf{Proof.}
%\strut
%\begin{enumerate}
%\item
Let us first verify assertion (i). Since $(q_{n},u_{n})\rightarrow(1,t\phi
_{1})$ in $\mathbb{R}\times W^{2,r}(\Omega)$ for some $t>0$, we have $\Phi
_{q}(1,t)=0$ by the implicit function theorem. By direct computations, we get
\begin{equation}
\Phi_{q}(q,t)=\int_{\Omega}a\left(  x\right)  \left[  (t\phi_{1}%
+w)^{q}\left\{  \log(t\phi_{1}+w)+\frac{qw_{q}}{t\phi_{1}+w}\right\}
-w_{q}\right]  \phi_{1}.\label{Pq}%
\end{equation}
Putting $q=1$ and using that $w(1,t)=0$, we find that
\begin{align}
\Phi_{q}(1,t) &  =\int_{\Omega}a\left(  x\right)  \left[  (t\phi_{1})\left\{
\log(t\phi_{1})+\frac{w_{q}(1,t)}{t\phi_{1}}\right\}  -w_{q}(1,t)\right]
\phi_{1}\nonumber\\
&  =t\int_{\Omega}a\left(  x\right)  \phi_{1}^{2}\log(t\phi_{1})\label{Pq1t}\\
&  =t\left\{  (\log t)\int_{\Omega}a\left(  x\right)  \phi_{1}^{2}%
+\int_{\Omega}a\left(  x\right)  \phi_{1}^{2}\log\phi_{1}\right\}  .\nonumber
\end{align}
Thus
\[
t=t_{\ast}:=\exp\left[  -\frac{\int_{\Omega}a\left(  x\right)  \phi_{1}%
^{2}\log\phi_{1}}{\int_{\Omega}a\left(  x\right)  \phi_{1}^{2}}\right]  ,
\]
as claimed in assertion (i). \newline

%Summing up, we have verified that if the unique positive solution $u_{q}$ of $(P_{a,q})$ converges to $t\phi_{1}$ for some $t>0$, then we have $t=t_{\ast}$.

Next, we verify assertion (ii). To this end, we use the fact that the map
$(q,t)\mapsto N(q,t)=t^{q}$ is analytic around $(q,t)=(1,t_{\ast})$, and apply
the implicit function theorem. We consider partial derivatives of $\Phi$, and
check that $\frac{\partial^{k}\Phi}{\partial t^{k}}(1,t)=0$ and $\Phi
_{qt}(1,t_{\ast})>0$. In fact, the case $k=1$ is straightforward since
$\Gamma_{1}$ is a trivial line of solutions of $(P_{a,q})$. Moreover, for
$k\geq2$, we have that
\[
\frac{\partial^{k}\Phi}{\partial t^{k}}(q,t)=(q-1)\Phi_{k}(q,t)+\int_{\Omega
}a(x)\left\{  q(t\phi_{1}+w)^{q-1}\frac{\partial^{k}w}{\partial t^{k}%
}(q,t)-\frac{\partial^{k}w}{\partial t^{k}}(q,t)\right\}  \phi_{1}%
\]
for some continuous function $\Phi_{k}$ of $(q,t)$ at $(1,t)$, so that
$\frac{\partial^{k}\Phi}{\partial t^{k}}(1,t)=0$ for all $k\in\mathbb{N}$ and
$t>0$. Since $(q,t)\mapsto t^{q} = \exp[q\log t]$ is analytic at $(q,t) =
(1,t)$, for any $t>0$, a regularity result for the implicit function theorem (see e.g.
\cite{Zei86}) ensures that so is $w(q,t)$ around $(1,t_{\ast})$, and thus so
is $\Phi(q,t)$. Combining this result with the fact that $\frac{\partial
^{k}\Phi}{\partial t^{k}}(1,t)=0$ for all $k\in\mathbb{N}$, we deduce that
$\Phi(q,t)$ is given around $(1,t_{\ast})$ by
\begin{align*}
&  \Phi(q,t)=(q-1)\hat{\Phi}(q,t),\quad\mbox{where}\\
&  \hat{\Phi}(q,t)=\frac{1}{2}\Phi_{qq}(1,t_{\ast})(q-1)+\Phi_{qt}(1,t_{\ast
})(t-t_{\ast})\\
&  \qquad\qquad+\mbox{ higher order terms w.r.t. }(q-1)\mbox{ and }(t-t_{\ast
}).
\end{align*}
Therefore, applying the implicit function theorem to $\hat{\Phi}(q,t)$ at
$(1,t_{\ast})$, we infer that the set $\Phi(q,t)=0$ around $(1,t_{\ast})$ is
given completely by
\[
q=1,\quad\mbox{and} \quad t=t(q)\quad\mbox{with} \quad t(1)=t_{\ast},
\]
provided that $\Phi_{qt}(1,t_{\ast})\not =0$. The desired conclusion follows.

Finally, we check that $\Phi_{qt}(1,t_{\ast})>0$: by a direct computation from
\eqref{Pq}, we observe that
\begin{align*}
\Phi_{qt} &  =\int_{\Omega}a\left(  x\right)  \left[  q(t\phi_{1}%
+w)^{q-1}(\phi_{1}+w_{t})\left\{  \log(t\phi_{1}+w)+\frac{qw_{q}}{t\phi_{1}%
+w}\right\}  \right.  \\
&  \quad\left.  +(t\phi_{1}+w)^{q}\left\{  \frac{\phi_{1}+w_{t}}{t\phi_{1}%
+w}+\frac{qw_{qt}(t\phi_{1}+w)-qw_{q}(\phi_{1}+w_{t})}{(t\phi_{1}+w)^{2}%
}\right\}  -w_{qt}\right]  \phi_{1}.
\end{align*}
	% \marginpar{\textit{modified}} 
Letting $q=1$, it follows that
\begin{equation}
\Phi_{qt}(1,t)=\int_{\Omega}a\left(  x\right)  [(\phi_{1}+w_{t}(1,t))\log
(t\phi_{1})+(\phi_{1}+w_{t}(1,t))]\phi_{1}.\label{Pqt}%
\end{equation}
We differentiate \eqref{beq01} with respect to $t$, and we obtain that
\[
-\Delta w_{t}-a\left(  x\right)  w_{t}=Q[a\left(  x\right)  \{q(t\phi
_{1}+w)^{q-1}(\phi_{1}+w_{t})-(\phi_{1}+w_{t})\}].
\]
Letting $q=1$ again, we deduce that
\[
-\Delta w_{t}(1,t)-a\left(  x\right)  w_{t}(1,t)=0,\ \mbox{ and }w_{t}(1,t)\in
X_{2}.
\]
Hence, $w_{t}(1,t)=0$, and thus, it follows from \eqref{Pqt} that
\begin{align*}
\Phi_{qt}(1,t) &  =\int_{\Omega}a\left(  x\right)  [\phi_{1}\log(t\phi
_{1})+\phi_{1}]\phi_{1}\\
&  =\int_{\Omega}a\left(  x\right)  \phi^{2}\log(t\phi_{1})+\int_{\Omega
}a\left(  x\right)  \phi_{1}^{2}%
\end{align*}
When $t=t_{\ast}$, we know that $\int_{\Omega}a\left(  x\right)  \phi^{2}%
\log(t_{\ast}\phi_{1})=0$ from \eqref{Pq1t}, so that
\[
\Phi_{qt}(1,t_{\ast})=\int_{\Omega}a\left(  x\right)  \phi_{1}^{2}>0,
\]
as desired. $\blacksquare$\newline

Next, as we did for $q<1$ close to $1$, we show that the Lyapunov-Schmidt
reduction is useful for the case $q>0$ close to $0$. Indeed, we
%\marginpar{\textit{modified}}
exhibit how to construct $a$ such that $(P_{a,q})$ possesses a solution in
$\mathcal{P}^{\circ}$ for $q>0$ arbitrarily close to $0$. Consider the problem%
\[
\left\{
\begin{array}
[c]{lll}%
-\Delta u=a(x) & \mathrm{in} & \Omega,\\
\frac{\partial u}{\partial\nu}=0 & \mathrm{on} & \partial\Omega.
\end{array}
\right.  \leqno{(P_{a,0})}
\]
It is easy to check that $(P_{a,0})$ has a solution if and only if
$\int_{\Omega}a=0$, in which case all solutions are of the form $u+c$, where
$c$ is any constant and $u$ is a particular solution.

Assume now that $a\not \equiv 0$ and $\int_{\Omega}a=0$ (in particular, $a$
changes sign). We set $X:=\{u\in W^{2,r}(\Omega):\frac{\partial u}{\partial
\nu}=0\mbox{ on }\partial\Omega\}$, and write $X=\langle1\rangle+X_{2}$, where
$\langle1\rangle$ is the set of constant functions and $X_{2}:=\{w\in
W^{2,r}(\Omega):\int_{\Omega}w=0\}$. Let $w_{0}$ be the unique solution of
$(P_{a,0})$ such that $w_{0}\in X_{2}$, and $t_{0}>0$ be such that
$u_{0}=t_{0}+w_{0}>0$ on $\overline{\Omega}$. Then $u_{0}\in\mathcal{P}%
^{\circ}$ solves $(P_{a,0})$.

Given $\varepsilon,\delta>0$ and $q\in(-\delta,\delta)$, we 
consider the following perturbation of $(P_{a,0})$: 
\[
\left\{
\begin{array}
[c]{lll}%
-\Delta u=(a(x)-\varepsilon)u^{q} & \mathrm{in} & \Omega,\\
\frac{\partial u}{\partial\nu}=0 & \mathrm{on} & \partial\Omega.
\end{array}
\right.  \leqno{(P_{a-\varepsilon,q})}
\]
Note that if $\varepsilon$ is sufficiently small, then $a-\varepsilon$ changes
sign, and $\int_{\Omega}(a(x)-\varepsilon)<0$. Note also that
$(P_{a-\varepsilon,q})$ admits $(q,\varepsilon,u)=(0,0,t_{0}+w_{0})$ as a
solution. Our aim is to look for positive solutions of $(P_{a-\varepsilon,q})$
in a neighborhood of $(0,0,t_{0}+w_{0})$. 
	% \marginpar{\it $Y_1$ is replaced by $Y_2$.}
Let $Y_{2}:=\{f\in L^{r}%
(\Omega):\int_{\Omega}f=0\}$ and $Q$ be the usual projection of $L^{r}%
(\Omega)$ to $Y_{2}$, given by $Q[f]:=f-\frac{1}{|\Omega|}\int_{\Omega}f$.
Following the Lyapunov-Schmidt approach already used in Theorem \ref{t1}, we
reduce $(P_{a-\varepsilon,q})$ to the following coupled equations
\begin{align}
&  Q[-\Delta u]=Q[(a-\varepsilon)u^{q}],\label{Q:eq}\\
&  (1-Q)[-\Delta u]=(1-Q)[(a-\varepsilon)u^{q}].\label{1-Q:eq}%
\end{align}
Associated with \eqref{Q:eq}, we define the mapping
\[
F:(-\delta_{0},\delta_{0})\times(-\varepsilon_{0},\varepsilon_{0})\times
(t_{0}-d_{0},t_{0}+d_{0})\times B_{\rho_{0}}(w_{0})\rightarrow Y_{2}%
\]
by
\[
F(q,\varepsilon,t,w):=-\Delta w-Q[(a-\varepsilon)(t+w)^{q}],
\]
where $B_{\rho_{0}}(w_{0})$ is the ball in $X_{2}$ with center $w_{0}$ and
radius $\rho_{0}$. We note that $F(0,0,t_{0},w_{0})=0$, since $Q[a]=a$. The
Fr\'{e}chet derivative $F_{w}(q,\varepsilon,t,w):X_{2}\rightarrow Y_{2}$ is
given by
\[
F_{w}(q,\varepsilon,t,w)\varphi=-\Delta\varphi-Q[(a-\varepsilon)q(t+w)^{q-1}%
\varphi].
\]
Taking $(q,\varepsilon,t,w)=(0,0,t_{0},w_{0})$, we see that $F_{w}%
(0,0,t_{0},w_{0})\varphi=-\Delta\varphi$, so that $F_{w}(0,0,t_{0},w_{0})$ is
bijective, and the implicit function theorem applies. Consequently, we have
\begin{align*}
&  F(q,\varepsilon,t,w)=0,\quad(q,\varepsilon,t,w)\simeq(0,0,t_{0},w_{0})\\
&  \Longleftrightarrow w=w(q,\varepsilon,t),\ \ (q,\varepsilon,t)\simeq
(0,0,t_{0})\ \mbox{ with }w(0,0,t_{0})=w_{0}.
\end{align*}
Using $w(q,\varepsilon,t)$, we derive from \eqref{1-Q:eq} the equation
\[
\Psi(q,\varepsilon,t):=\int_{\Omega}(a(x)-\varepsilon)(t+w(q,\varepsilon
,t))^{q}=0\quad\mbox{ in }\mathbb{R}^{3}.
\]
Note that $\Psi(0,0,t_{0})=0$. We prove now the following result:

\begin{proposition}
\label{prop:q0} Given $\varepsilon>0$ sufficiently small, there exists
$q_{\varepsilon}>0$ such that $u_{\varepsilon}=t_{0}+w(q_{\varepsilon
},\varepsilon,t_{0})>0$ on $\overline{\Omega}$, and $\Psi(q_{\varepsilon
},\varepsilon,t_{0})=0$. Moreover, $q_{\varepsilon}\rightarrow0^{+}$ as
$\varepsilon\rightarrow0^{+}$. Consequently, $u_{\varepsilon}\in
\mathcal{P}^{\circ}$ is a solution of $(P_{a-\varepsilon,q_{\varepsilon}})$.
\newline
\end{proposition}

\textbf{Proof.} We apply the implicit function theorem for $\Psi$ at
$(0,0,t_{0})$. We observe that
\[
\frac{\partial\Psi}{\partial q}=\int_{\Omega}(a(x)-\varepsilon)(t+w)^{q}%
\left\{  \log(t+w)+\frac{q}{t+w}\frac{\partial w}{\partial q}\right\},
\]
and therefore
\[
\frac{\partial\Psi}{\partial q}(0,0,t_{0})=\int_{\Omega}a(x)\log u_{0}.
\]
Since $u_{0}$ is a solution in $\mathcal{P}^{\circ}$ of $(P_{a,0})$, we see
that
\[
0<\int_{\Omega}\frac{|\nabla u_{0}|^{2}}{u_{0}}=\int_{\Omega}-\Delta u_{0}\log
u_{0}=\int_{\Omega}a(x)\log u_{0},
\]
which implies that $\frac{\partial\Psi}{\partial q}(0,0,t_{0})>0$. Hence, the
implicit function theorem ensures that
\[
\Psi(q,\varepsilon,t)=0,\ \ (q,\varepsilon,t)\simeq(0,0,t_{0})\Leftrightarrow
q=q(\varepsilon,t),\ \ (\varepsilon,t)\simeq(0,t_{0}%
)\ \mbox{ with }\ q(0,t_{0})=0.
\]

Next we show that $\frac{\partial q}{\partial\varepsilon}(0,t_{0})>0$. To this
end, we differentiate $\Psi$ with respect to $\varepsilon$, obtaining that
\[
\frac{\partial\Psi}{\partial\varepsilon}=\int_{\Omega}\left\{  -(t+w)^{q}%
+(a(x)-\varepsilon)q(t+w)^{q-1}\frac{\partial w}{\partial\varepsilon}\right\},
\]
and so 
\[
\frac{\partial\Psi}{\partial\varepsilon}(0,0,t_{0})=-|\Omega|<0.
\]
It follows that
\[
\frac{\partial q}{\partial\varepsilon}(0,t_{0})=-\frac{\frac{\partial\Psi
}{\partial\varepsilon}(0,0,t_{0})}{\frac{\partial\Psi}{\partial q}(0,0,t_{0}%
)}=\frac{|\Omega|}{\int a(x)\log u_{0}}>0.
\]
Using the mean value theorem, we deduce that for $\varepsilon>0$ small
enough,
\[
q(\varepsilon,t_{0})=q(0,t_{0})+\frac{\partial q}{\partial\varepsilon}%
(\theta\varepsilon,t_{0})\varepsilon>0,
\]
for some $0<\theta<1$. Thus
\[
q(\varepsilon,t_{0})\rightarrow q(0,t_{0})=0\quad\mbox{ as }\varepsilon
\rightarrow0^{+},
\]
as desired. $\blacksquare$ \newline

\begin{remark}
Let us analyze the asymptotic behavior of nontrivial nonnegative solutions of
$(P_{a,q})$ as $q\rightarrow0^{+}$ under $(H_{0})$ and $(H_{1})$. From
Proposition \ref{t2.1}, we know that bifurcation from zero or from infinity
does not occur as $q\rightarrow0^{+}$. It is thus natural to investigate the
limit of a sequence $\{u_{n}\}$ of nontrivial nonnegative solutions of
$(P_{a,q_{n}})$ with $q_{n}\rightarrow0^{+}$. Since $\{u_{n}\}$ is bounded in
$L^{\infty}(\Omega)$, it follows, by elliptic regularity, that up to a
subsequence, $u_{n}\rightarrow u_{0}$ in $\mathcal{C}^{1}(\overline{\Omega})$
with $u_{0}\not \equiv 0$. 
We point out that $u_{0}$ must vanish in a nonempty subset of $\Omega$ 
	% \marginpar{\it modified} 
with \textit{positive measure} 
(in other words, $u_{0}$ has a nonempty dead core). 
Indeed, if $u_{0}>0$ a.e. in $\Omega$, then, passing to the limit, 
we have that
\[
\int_{\Omega}\nabla u_{0}\nabla v=\int_{\Omega}a(x)v,\quad\forall
v\in\mathcal{C}^{1}(\overline{\Omega}),
\]
i.e. $u_{0}$ is a positive solution of $(P_{a,0})$. Integrating this equation,
we deduce that $\int_{\Omega}a=0$, which is a contradiction.
\end{remark}

\bigskip

\subsection{Stability properties}

We conclude this section discussing the stability of the bifurcating positive
solutions provided by Theorem \ref{t1} (ii).

\begin{proposition}
\label{t2} Assume $(H_{0})$. If $\mu_{1}(a)=1$, then the bifurcating positive
solution $u(q)=t(q)\phi_{1}+w(q,t(q))$ given by Theorem \ref{t1} (ii) is
asymptotically stable (respect.\ unstable) for $q<1$ (respect.\ $q>1$).
\end{proposition}

\textbf{Proof.} Consider
\begin{equation}
-\Delta\varphi_{1}(q)=qa(x)u(q)^{q-1}\varphi_{1}(q)+\gamma_{1}(q)\varphi
_{1}(q),\label{e2}%
\end{equation}
where $\gamma_{1}(q):=\gamma_{1}(q,u(q))$, and $\varphi_{1}(q)$ is a positive
eigenfunction associated to $\gamma_{1}(q)$. We see that $\gamma_{1}(1)=0$ and
$\varphi_{1}(1)=\phi_{1}$. To analyse $\gamma_{1}(q)$ for $q\neq1$, we
differentiate \eqref{e2} with respect to $q$, to obtain that
\begin{align*}
-\Delta\varphi_{1}^{\prime} &  =a\left(  x\right)  u^{q-1}\varphi
_{1}+qa(x)u^{q-1}\left(  \log u+(q-1)\frac{u^{\prime}}{u}\right)  \varphi
_{1}+qa(x)u^{q-1}\varphi_{1}^{\prime}\\
&  \quad+\gamma_{1}^{\prime}\varphi+\gamma_{1}\varphi_{1}^{\prime}.
\end{align*}
Letting $q=1$ here, it follows that
\[
A\varphi_{1}^{\prime}(1)=\gamma_{1}^{\prime}(1)\phi_{1}+a(x)\left\{  \phi
_{1}+\phi_{1}\log(t_{\ast}\phi_{1})\right\}  ,
\]
	% \marginpar{\textit{modified}} 
and thus, by the divergence theorem,%
\[
0=\int_{\Omega}A\varphi_{1}^{\prime}(1)\phi_{1}-\varphi_{1}^{\prime}%
(1)A\phi_{1}=\gamma_{1}^{\prime}(1)+\int_{\Omega}a(x)\left(  \phi_{1}+\phi
_{1}\log(t_{\ast}\phi_{1})\right)  \phi_{1}.
\]
Since $\int_{\Omega}a(x)\phi_{1}^{2}\log(t_{\ast}\phi_{1})=0$, we obtain that
\[
\gamma_{1}^{\prime}(1)=-\int_{\Omega}a(x)\phi_{1}^{2}<0.
\]
The desired conclusion follows from the fact that $\gamma_{1}(1)=0$.
$\blacksquare$

\begin{remark}
\label{rem:stab} 
\strut
\begin{enumerate}
\item The stability result of Proposition \ref{t2} also follows from \cite[Theorem 1]{BH}. Even though this result assumes $a$ to be smooth and the nonlinearity to be $C^2$ at $0$, one may easily see that under our assumptions it also applies to solutions of $(P_{a,q})$ in $\mathcal{P}^{\circ}$. More generally, it shows that any such solution is asymptotically stable for \textit{every} $0<q<1$.\\

\item When $q>1$, we can deduce (by a well known approach) that
\textit{every} solution $u\in\mathcal{P}^{\circ}$ of $(P_{a,q})$ is unstable.
Indeed, linearizing $(P_{a,q})$ at $u$ we obtain $-\Delta\varphi
=qa(x)u^{q-1}\varphi+\gamma\varphi$. The divergence theorem yields that
\[
0>\int_{\Omega}\frac{u}{\varphi_{1}}\sum_{j}\frac{\partial}{\partial x_{j}%
}\left(  \varphi_{1}^{2}\frac{\partial}{\partial x_{j}}\left(  \frac
{u}{\varphi_{1}}\right)  \right)  =(q-1)\int_{\Omega}a(x)u^{q+1}+\gamma
_{1}\int_{\Omega}u^{2},
\]
where $\gamma_{1}=\gamma_{1}(q)$ and $\varphi_{1}=\varphi_{1}(q)$. Hence, we
obtain
\[
\gamma_{1}<\frac{(1-q)\int_{\Omega}a(x)u^{q+1}}{\int_{\Omega}u^{2}}<0.
\]

\end{enumerate}
\end{remark}

\bigskip

\section{Proofs of Theorems \ref{c1}, \ref{cc} and \ref{rad2}}

\textbf{Proof of Theorem \ref{c1}:} Let us first observe that by Theorem
\ref{t1}, there exists $q_{0}=q_{0}(a)<1$ such that $(P_{a,q})$ has a solution
$u_{q}\in\mathcal{P}^{\circ}$ for $q_{0}<q<1$. Moreover,
	% \marginpar{\textit{modified}}
the proof of \cite[Lemma 3.1]{BPT} can be adapted
to our setting, so that $(P_{a,q})$ has no other positive solution for
$q_{0}<q<1$. We consider now the asymptotic behavior of $u_{q}$ as
$q\rightarrow1^{-}$. Assertion (i) is a direct consequence of Theorem \ref{t1}
(ii) and elliptic regularity.

Assume now that $\mu_{1}=\mu_{1}(a)\neq1$ and set $v:=\mu_{1}^{\frac{1}{1-q}%
}u$. Note that if $u$ solves $(P_{a,q})$ then $v$ solves $(P_{\tilde{a},q})$,
where $\tilde{a}:=\mu_{1}a$. Indeed,
\[
-\Delta v=\mu_{1}^{\frac{1}{1-q}}a\left(  x\right)  u^{q}=\mu_{1}a\left(
x\right)  v^{q}=\tilde{a}\left(  x\right)  v^{q}.
\]
Moreover, we easily see that $\mu_{1}(\tilde{a})=1$. By item (i), we get a
positive solution $v_{q}$ of $(P_{\tilde{a},q})$ such that $v_{q}\rightarrow
t_{\ast}(\tilde{a})\phi_{1}(\tilde{a})$, where $\phi_{1}(\tilde{a})$ is a
positive eigenfunction of $(E_{1,\tilde{a}})$, which is nothing but
$(E_{\mu_{1},a})$, i.e. $\phi_{1}(\tilde{a})=\phi_{1}(a)$ and $t_{\ast}%
(\tilde{a})=t_{\ast}(a)$ . In this way, we obtain a positive solution
$u_{q}=\mu_{1}^{\frac{1}{q-1}}v_{q}$ of $(P_{a,q})$ for $q$ close to $1$. In
particular, we see that if $\mu_{1}>1$ then $\mu_{1}^{\frac{1}{q-1}%
}\rightarrow0$, so that $u_{q}\rightarrow0$ in $\mathcal{C}^{1}(\overline
{\Omega})$ as $q\rightarrow1^{-}$. On the other hand, if $\mu_{1}<1$, then
$\mu_{1}^{\frac{1}{q-1}}\rightarrow\infty$, so that $\displaystyle\min
_{\overline{\Omega}}u_{q}\rightarrow\infty$ 
	% \marginpar{\textit{modified}} 
when $q\rightarrow1^{-}$.

Finally, the asymptotic stability of $u_{q}$ is a direct consequence of
Proposition \ref{t2}. $\blacksquare$\newline

When proving Theorems \ref{cc}, \ref{rad2} and \ref{ti}, we shall repeatedly
use the following remark:

\begin{remark}
\label{rema} \strut

\begin{enumerate}
\item[(i)] Since $(P_{a,q})$ is homogeneous, we see that $(P_{a,q})$ has a
nonnegative (respect. positive) solution if and only if, for any $\sigma>0$
fixed, $(P_{\sigma a,q})$ has a nonnegative (respect. positive) solution.

\item[(ii)] Lemma 2.4 in \cite{BPT} (which is proved using Proposition 2.1
therein) gives the existence of arbitrarily large supersolutions of
$(P_{a,q})$ provided that $\int_{\Omega}a<0$. Although it is assumed that $a$
is H\"{o}lder continuous in \cite{BPT}, one can see that Lemma 2.4 and
Proposition 2.1 still hold (with the same proof) if $a\in L^{\infty}(\Omega)$.
\end{enumerate}
\end{remark}

\medskip

\textbf{Proof of Theorem \ref{cc}:} We proceed in several steps. By Remark
\ref{rema}, it is enough to provide a positive (in $\Omega$) weak subsolution for
$(P_{b,q})$, where $b:=\gamma a$ and $\gamma:=1/\left(  1-q\right)  $.
%\marginpar{\textit{added}}
Observe that $\gamma q=\gamma-1$. We note also that, since $\int_{B_{R}}a<0$,
it holds that $R_{0}<R$. Let us first define%
\begin{align*}
C  &  :=\frac{1-q}{1+q},\\
w\left(  r\right)   &  :=C\int_{r}^{R}\frac{1}{t^{N-1}}\int_{t}^{R}%
a^{-}\left(  y\right)  y^{N-1}\,dy\,dt:=C\phi\left(  r\right)  ,\quad
r\in\left[  R_{0},R\right]  .
\end{align*}

Then, $w\left(  R\right)  =w^{\prime}\left(  R\right)  =0$ and $w\left(
r\right)  >0$ for all $r\in\left[  R_{0},R\right)  $. Also, a few computations
show that
\begin{equation}
\phi^{\prime\prime}+\frac{N-1}{r}\phi^{\prime}=a^{-}\left(  r\right)  .
\label{pola}%
\end{equation}
Let now $z\left(  r\right)  :=w^{\gamma}\left(  r\right)  $. We claim that
\begin{equation}
-\Delta z\leq\gamma a\left(  x\right)  z^{q}\quad a.e.\text{ in } A_{R_{0}%
,R}. \label{ani}%
\end{equation}
Indeed, since $z$ is radial, there holds
\begin{align*}
\Delta z  &  =z^{\prime\prime}+\frac{N-1}{r}z^{\prime}\\
&  =\gamma\left(  \gamma-1\right)  w^{\gamma-2}\left(  w^{\prime}\right)
^{2}+\gamma w^{\gamma-1}w^{\prime\prime}+\frac{\gamma\left(  N-1\right)  }%
{r}w^{\gamma-1}w^{\prime},
\end{align*}
and also
\[
-\gamma a\left(  r\right)  z^{q}\leq\gamma a^{-}\left(  r\right)  w^{\gamma
q}=\gamma a^{-}\left(  r\right)  w^{\gamma-1}.
\]
Thus, in order to prove the claim it is enough to verify that
\[
\left(  \gamma-1\right)  \frac{\left(  w^{\prime}\right)  ^{2}}{w}%
+w^{\prime\prime}+\frac{N-1}{r}w^{\prime}\geq a^{-}\left(  r\right)  .
\]
Now, taking into account (\ref{pola}) and that $w=C\phi$, the above inequality
is equivalent to
\begin{equation}
F\left(  r\right)  :=\left(  \gamma-1\right)  \left(  \phi^{\prime}\right)
^{2}\geq\left(  \frac{1}{C}-1\right)  a^{-}\left(  r\right)  \phi:=G\left(
r\right)  . \label{fg}%
\end{equation}
We observe next that $F\left(  R\right)  =G\left(  R\right)  =0$ and
$F^{\prime}\left(  r\right)  \leq0$ for all $r\in\left[  R_{0},R\right]  $
(recall that $\phi^{\prime}\leq0$). So, in order to check (\ref{fg}) it
suffices to see that $F^{\prime}\left(  r\right)  \leq G^{\prime}\left(
r\right)  $ for such $r$. Now,%
\begin{align*}
F^{\prime}\left(  r\right)   &  =2\left(  \gamma-1\right)  \phi^{\prime}%
\phi^{\prime\prime},\\
G^{\prime}\left(  r\right)   &  =\left(  \frac{1}{C}-1\right)  \left(  \left(
a^{-}\left(  r\right)  \right)  ^{\prime}\phi+a^{-}\left(  r\right)
\phi^{\prime}\right)  \geq\left(  \frac{1}{C}-1\right)  a^{-}\left(  r\right)
\phi^{\prime},
\end{align*}
where we used the fact that $a$ is differentiable and nonincreasing in
$A_{R_{0},R}$. Therefore, $F^{\prime}\left(  r\right)  \leq G^{\prime}\left(
r\right)  $ provided that%
\[
2\left(  \gamma-1\right)  \phi^{\prime}\phi^{\prime\prime}\leq\left(  \frac
{1}{C}-1\right)  a^{-}\left(  r\right)  \phi^{\prime},
\]
i.e.%
\begin{equation}
2\left(  \gamma-1\right)  \left(  -\frac{N-1}{r}\phi^{\prime}+a^{-}\left(
r\right)  \right)  \geq\left(  \frac{1}{C}-1\right)  a^{-}\left(  r\right)  .
\label{puf}%
\end{equation}
But (\ref{puf}) holds by our election of $C$. Indeed, since $\phi^{\prime}%
\leq0$, one only has to observe that $2\left(  \gamma-1\right)  =\frac{1}%
{C}-1.$

On the other side, let $v$ be a solution of%
\begin{equation}
\left\{
\begin{array}
[c]{lll}%
-\Delta v=\gamma a\left(  x\right)  v^{q} & \mathrm{in} & B_{R_{0}},\\
v=z\left(  R_{0}\right)  & \mathrm{on} & \partial B_{R_{0}}.
\end{array}
\right.  \label{veve}%
\end{equation}
Such $v$ can be easily constructed by the sub and supersolutions method, since
$a\geq0$ in $B_{R_{0}}$. Moreover, $v$ is radial. Indeed, this follows from
either the fact that the sub and supersolutions can be chosen radial, or
because the solution of (\ref{veve}) is unique (cf. \cite{DS}) and $v\left(
Sx\right)  $ is also a solution if $S$ is an isometry of $\mathbb{R}^{N}$.
Furthermore, it is also easy to check that $r\rightarrow v(r)$ is
nonincreasing in $(0,R_{0})$ because $a\geq0$ in $B_{R_{0}}$. Hence, by the
divergence theorem (as stated e.g. in \cite{cuesta}, p. 742),
\begin{align}
v^{\prime}\left(  R_{0}\right)  \omega_{N-1}R_{0}^{N-1}  &  =\int_{B_{R_{0}}%
}\Delta v=-\int_{B_{R_{0}}}\gamma av^{q}\nonumber\\
&  \leq-\gamma v^{q}\left(  R_{0}\right)  \int_{B_{R_{0}}}a=-\gamma w^{\gamma
q}\left(  R_{0}\right)  \int_{B_{R_{0}}}a. \label{uuu}%
\end{align}
On the other hand, recalling that $\gamma-1=\gamma q$, we obtain that
\[
z^{\prime}\left(  R_{0}\right)  =\gamma w^{\gamma-1}\left(  R_{0}\right)
w^{\prime}\left(  R_{0}\right)  =-\gamma w^{\gamma q}\left(  R_{0}\right)
\frac{C}{R_{0}^{N-1}}\int_{R_{0}}^{R}a^{-}\left(  y\right)  y^{N-1}\,dy
\]
and so
\begin{align}
z^{\prime}\left(  R_{0}\right)  \omega_{N-1}R_{0}^{N-1}  &  =-\gamma w^{\gamma
q}\left(  R_{0}\right)  C\omega_{N-1}\int_{R_{0}}^{R}a^{-}\left(  y\right)
y^{N-1}\,dy\nonumber\\
&  =-\gamma w^{\gamma q}\left(  R_{0}\right)  C\int_{A_{R_{0}}}a^{-}.
\label{dos}%
\end{align}
Next we observe that $v^{\prime}\left(  R_{0}\right)  \leq z^{\prime}\left(
R_{0}\right)  $. Indeed, taking into account \eqref{uuu}, \eqref{dos} and
the definition of $C$, we see that this is true by \eqref{inferno}.

To conclude the existence assertion, we define $u:=z$ in $\overline{A}%
_{R_{0},R}$ and $u:=v$ in $B_{R_{0}}$. Then $u\in H_{0}^{1}\left(
\Omega\right)  \cap L^{\infty}\left(  \Omega\right)  $, $u>0$ in $\Omega$ and
$\frac{\partial u}{\partial\nu}=0$ on $\partial\Omega$. Moreover, recalling
(\ref{ani}), (\ref{veve}) and that $v^{\prime}\left(  R_{0}\right)  \leq
z^{\prime}\left(  R_{0}\right)  $, the divergence theorem yields that
$u$ is a weak subsolution of $(P_{b,q})$.

Finally, the uniqueness assertion
%\marginpar{\textit{modified}}
is a consequence of \cite[Theorem 3.1]{BPT} (see Remark \ref{h}).
$\blacksquare$\newline\medskip

\textbf{Proof of Theorem \ref{rad2}:} Given $\varepsilon\in\left[  0,q\right)
$, we define%
\[
\gamma_{\varepsilon}:=\frac{1-\varepsilon}{1-q},\quad C_{\varepsilon}:=\left(
\frac{R_{0}^{2\varepsilon}}{2}\frac{\left\Vert a^{-}\right\Vert _{\infty}%
}{2\left(  \gamma_{\varepsilon}-1\right)  +N}+\varepsilon\right)  ^{\frac
{1}{1-\varepsilon}}.
\]
We note that both $\varepsilon\mapsto\gamma_{\varepsilon}$ and $\varepsilon
\mapsto C_{\varepsilon}$ are continuous, that
\begin{equation}
\varepsilon=\gamma_{\varepsilon}q-\left(  \gamma_{\varepsilon}-1\right)
,\quad\gamma_{\varepsilon}-1=\frac{q-\varepsilon}{1-q}>0, \label{ga}%
\end{equation}
and that%
\begin{equation}
2C_{0}=\frac{\left\Vert a^{-}\right\Vert _{\infty}}{2\left(  \gamma
_{0}-1\right)  +N}=\frac{1-q}{2q+N\left(  1-q\right)  }\left\Vert
a^{-}\right\Vert _{\infty}. \label{c0}%
\end{equation}

Given $r\in\left[  0,R_{0}\right]  $ and $\delta\geq0$, we set $u_{\delta
,\varepsilon}\left(  r\right)  :=C_{\varepsilon}r^{2}+\delta$. We now observe
that we can fix $\varepsilon>0$ small enough such that
\begin{equation}
u_{0,\varepsilon}^{\prime}\left(  R_{0}\right)  <u_{0,\varepsilon
}^{\varepsilon}\left(  R_{0}\right)  \frac{1}{R_{0}^{N-1}}\int_{R_{0}}%
^{R}a\left(  y\right)  y^{N-1}\,dy. \label{a}%
\end{equation}
Indeed, $u_{0,\varepsilon}^{\prime}\left(  R_{0}\right)  =2C_{\varepsilon
}R_{0}$ and $u_{0,\varepsilon}^{\varepsilon}\left(  R_{0}\right)  =\left(
C_{\varepsilon}R_{0}^{2}\right)  ^{\varepsilon}$, and so (\ref{a}) holds if
and only if
\begin{equation}
2C_{\varepsilon}^{1-\varepsilon}R_{0}^{N}<R_{0}^{2\varepsilon}\int_{R_{0}}%
^{R}a\left(  y\right)  y^{N-1}\,dy. \label{b}%
\end{equation}
Now, by (\ref{sipi}) and (\ref{c0}),
\[
2C_{0}R_{0}^{N}=\frac{1-q}{2q+N\left(  1-q\right)  }R_{0}^{N}\left\Vert
a^{-}\right\Vert _{\infty}<\frac{\int_{A_{R_{0},R}}a}{\omega_{N-1}}%
=\int_{R_{0}}^{R}a\left(  y\right)  y^{N-1}\,dy.
\]
Thus, (\ref{b}) (and consequently (\ref{a})) holds for $\varepsilon>0$
sufficiently small.

Next, we note that, by definition,
\[
2C_{\varepsilon}^{1-\varepsilon}>\frac{\left\Vert a^{-}\right\Vert _{\infty
}R_{0}^{2\varepsilon}}{2\left(  \gamma_{\varepsilon}-1\right)  +N}.
\]
Therefore, for all $r\in\left[  0,R_{0}\right]  $,%
\begin{equation}
4C_{\varepsilon}\left(  \gamma_{\varepsilon}-1\right)  +2NC_{\varepsilon
}-\left\Vert a^{-}\right\Vert _{\infty}\left(  C_{\varepsilon}r^{2}\right)
^{\varepsilon}>0. \label{asa}%
\end{equation}
In view of this inequality, we may fix $\delta>0$ such that%
\begin{equation}
\frac{4C_{\varepsilon}^{2}r^{2}\left(  \gamma_{\varepsilon}-1\right)
}{C_{\varepsilon}r^{2}+\delta}+2NC_{\varepsilon}-\left\Vert a^{-}\right\Vert
_{\infty}\left(  C_{\varepsilon}r^{2}+\delta\right)  ^{\varepsilon}%
>0\quad\forall r\in\left[  0,R_{0}\right]  . \label{casa}%
\end{equation}
Indeed, we pick first any $\delta_{0}>0$ small enough such that, for all
$\delta\in\left(  0,\delta_{0}\right]  $,
\[
2NC_{\varepsilon}>\left\Vert a^{-}\right\Vert _{\infty}\delta^{\varepsilon}.
\]
Then there exists $r_{0}=r_{0}\left(  \delta_{0}\right)  >0$ such that
\[
2NC_{\varepsilon}>\left\Vert a^{-}\right\Vert _{\infty}\left(  C_{\varepsilon
}r^{2}+\delta\right)  ^{\varepsilon}\quad\forall r\in\left[  0,r_{0}\right],
\]
and thus (\ref{casa}) clearly holds for all $r\in\left[  0,r_{0}\right]  $.
Suppose now that $r\in\left[  r_{0},R_{0}\right]  $. Then, from (\ref{asa}) we
derive that
\begin{equation}
\frac{4C_{\varepsilon}^{2}r^{2}\left(  \gamma_{\varepsilon}-1\right)
}{C_{\varepsilon}r^{2}}+2NC_{\varepsilon}-\left\Vert a^{-}\right\Vert
_{\infty}\left(  C_{\varepsilon}r^{2}\right)  ^{\varepsilon}>0\quad\forall
r\in\left[  r_{0},R_{0}\right]  . \label{as}%
\end{equation}
Now, since by Dini's theorem the left-hand side of (\ref{casa}) converges to
the left-hand side of (\ref{as}) uniformly in $r\in\left[  r_{0},R_{0}\right]
$ as $\delta\rightarrow0^{+}$, then, decreasing $\delta_{0}$ if necessary, we
also see that (\ref{casa}) holds for all $r\in\left[  r_{0},R_{0}\right]  $.

Finally, since $u_{0,\varepsilon}^{\prime}=u_{\delta,\varepsilon}^{\prime}$
and $u_{0,\varepsilon}<u_{\delta,\varepsilon}$, recalling (\ref{a}), we get
that
\begin{equation}
u_{\delta,\varepsilon}^{\prime}\left(  R_{0}\right)  <u_{\delta,\varepsilon
}^{\varepsilon}\left(  R_{0}\right)  \frac{1}{R_{0}^{N-1}}\int_{R_{0}}%
^{R}a\left(  y\right)  y^{N-1}\,dy. \label{ad}%
\end{equation}
We fix for the rest of the proof $\varepsilon,\delta>0$ such that (\ref{casa})
and (\ref{ad}) hold.

Let $z_{\delta,\varepsilon}\left(  r\right)  :=u_{\delta,\varepsilon}%
^{\gamma_{\varepsilon}}\left(  r\right)  $. Let us show that%
\begin{equation}
\Delta z_{\delta,\varepsilon}\geq\gamma_{\varepsilon}\left\Vert a^{-}%
\right\Vert _{\infty}z_{\delta,\varepsilon}^{q}\quad\text{in }B_{R_{0}}.
\label{laaa}%
\end{equation}
Note that (\ref{laaa}) implies that $-\Delta z_{\delta,\varepsilon}\leq
\gamma_{\varepsilon}a\left(  x\right)  z_{\delta,\varepsilon}^{q}$ $a.e.$ in
$B_{R_{0}}$. We compute
%\marginpar{\textit{modified the above formula}} %
\begin{align*}
\Delta z_{\delta,\varepsilon}  &  =\gamma_{\varepsilon}\left(  \left(
\gamma_{\varepsilon}-1\right)  u_{\delta,\varepsilon}^{\gamma_{\varepsilon}%
-2}\left\vert \nabla u_{\delta,\varepsilon}\right\vert ^{2}+u_{\delta
,\varepsilon}^{\gamma-1}\Delta u_{\delta,\varepsilon}\right) \\
&  =\gamma_{\varepsilon}\left(  4C_{\varepsilon}^{2}r^{2}\left(
\gamma_{\varepsilon}-1\right)  u_{\delta,\varepsilon}^{\gamma_{\varepsilon}%
-2}+2NC_{\varepsilon}u_{\delta,\varepsilon}^{\gamma_{\varepsilon}-1}\right)  .
\end{align*}
Thus, in order to prove (\ref{laaa}) it is enough to see that%
\[
4C_{\varepsilon}^{2}r^{2}\left(  \gamma_{\varepsilon}-1\right)  u_{\delta
,\varepsilon}^{\gamma_{\varepsilon}-2}+2NC_{\varepsilon}u_{\delta,\varepsilon
}^{\gamma_{\varepsilon}-1}\geq\left\Vert a^{-}\right\Vert _{\infty}%
u_{\delta,\varepsilon}^{\gamma_{\varepsilon}q}.
\]
Furthermore, since $\varepsilon=\gamma_{\varepsilon}q-\left(  \gamma
_{\varepsilon}-1\right)  $ (recall (\ref{ga})), this is equivalent to
\[
\frac{4C_{\varepsilon}^{2}r^{2}\left(  \gamma_{\varepsilon}-1\right)
}{u_{\delta,\varepsilon}}+2NC_{\varepsilon}\geq\left\Vert a^{-}\right\Vert
_{\infty}u_{\delta,\varepsilon}^{\varepsilon}.
\]
But taking into account the definition of $u_{\delta,\varepsilon}$, we see
that the above inequality holds thanks to (\ref{casa}).

On the other side, let us define%
\begin{align*}
\phi\left(  r\right)   &  :=\int_{r}^{R}\frac{1}{t^{N-1}}\int_{t}^{R}a\left(
y\right)  y^{N-1}\,dy\,dt,\quad r\in\left[  R_{0},R\right]  ,\\
K  &  :=\frac{\gamma_{\varepsilon}}{\gamma_{0}}\phi\left(  R_{0}\right)
+\left[  u_{\delta,\varepsilon}\left(  R_{0}\right)  \right]  ^{\gamma
_{\varepsilon}/\gamma_{0}},\\
w\left(  r\right)   &  :=K-\frac{\gamma_{\varepsilon}}{\gamma_{0}}\phi\left(
r\right)  ,\quad\text{and}\\
v\left(  r\right)   &  :=w^{\gamma_{0}}\left(  r\right)  .
\end{align*}
Note that $v\left(  R_{0}\right)  =z_{\delta,\varepsilon}\left(  R_{0}\right)
$, and observe also that
\[
w^{\prime}\left(  r\right)  =\frac{\gamma_{\varepsilon}}{\gamma_{0}}\frac
{1}{r^{N-1}}\int_{r}^{R}a\left(  y\right)  y^{N-1}\,dy,
\]
and hence $w^{\prime}\left(  R\right)  =0$ (and $v^{\prime}\left(  R\right)
=0$). We also infer that $w\left(  r\right)  >0$ for $r\in\left[
R_{0},R\right]  $ since $w\left(  R_{0}\right)  >0$ and $w$ is increasing.
Moreover,
\begin{equation}
w^{\prime\prime}+\frac{N-1}{r}w^{\prime}=-\frac{\gamma_{\varepsilon}}%
{\gamma_{0}}a\left(  r\right)  . \label{polo}%
\end{equation}

We prove now that
\begin{equation}
-\Delta v\leq\gamma_{\varepsilon}a\left(  x\right)  v^{q}\quad a.e.\text{ in 
} A_{R_{0},R}. \label{ann}%
\end{equation}
Indeed, since $v$ is radial, (\ref{ann}) is equivalent to
\begin{align*}
-\Delta v  &  =-v^{\prime\prime}-\frac{N-1}{r}v^{\prime}\\
&  =-\gamma_{0}\left(  \left(  \gamma_{0}-1\right)  w^{\gamma_{0}-2}\left(
w^{\prime}\right)  ^{2}+w^{\gamma_{0}-1}w^{\prime\prime}+\frac{\left(
N-1\right)  }{r}w^{\gamma_{0}-1}w^{\prime}\right) \\
&  \leq\gamma_{\varepsilon}a\left(  r\right)  w^{\gamma_{0}-1}=\gamma
_{\varepsilon}a\left(  r\right)  w^{\gamma_{0}q}=\gamma_{\varepsilon}a\left(
r\right)  v^{q},
\end{align*}
and the above inequality clearly holds by (\ref{polo}).

We next verify that%
\begin{equation}
z_{\delta,\varepsilon}^{\prime}\left(  R_{0}\right)  \leq v^{\prime}\left(
R_{0}\right)  . \label{fin}%
\end{equation}
We have that
\begin{align*}
z_{\delta,\varepsilon}^{\prime}\left(  R_{0}\right)   &  =\gamma_{\varepsilon
}u_{\delta,\varepsilon}^{\gamma_{\varepsilon}-1}\left(  R_{0}\right)
u_{\delta,\varepsilon}^{\prime}\left(  R_{0}\right)  ,\\
v^{\prime}\left(  R_{0}\right)   &  =\gamma_{0}w^{\gamma_{0}-1}\left(
R_{0}\right)  w^{\prime}\left(  R_{0}\right)  ,
\end{align*}
and so it suffices to check that%
\begin{equation}
\gamma_{\varepsilon}u_{\delta,\varepsilon}^{\gamma_{\varepsilon}-1}\left(
R_{0}\right)  u_{\delta,\varepsilon}^{\prime}\left(  R_{0}\right)  \leq
\gamma_{0}w^{\gamma_{0}-1}\left(  R_{0}\right)  w^{\prime}\left(
R_{0}\right)  . \label{cal}%
\end{equation}
We observe now that, by definition, $w^{\gamma_{0}}\left(  R_{0}\right)
=u_{\delta,\varepsilon}^{\gamma_{\varepsilon}}\left(  R_{0}\right)  $, and
hence
\[
w^{\gamma_{0}-1}\left(  R_{0}\right)  =\left[  u_{\delta,\varepsilon}\left(
R_{0}\right)  \right]  ^{\frac{\gamma_{\varepsilon}\left(  \gamma
_{0}-1\right)  }{\gamma_{0}}}.
\]
Therefore, (\ref{cal}) can be written as
\[
\gamma_{\varepsilon}u_{\delta,\varepsilon}^{\prime}\left(  R_{0}\right)
\leq\gamma_{0}\left[  u_{\delta,\varepsilon}\left(  R_{0}\right)  \right]
^{\frac{\gamma_{\varepsilon}\left(  \gamma_{0}-1\right)  }{\gamma_{0}}-\left(
\gamma_{\varepsilon}-1\right)  }w^{\prime}\left(  R_{0}\right)  .
\]
Now, $\gamma_{0}-1=q/\left(  1-q\right)  $ and so, recalling the first
equality in (\ref{ga}), we see that%
\[
\frac{\gamma_{\varepsilon}\left(  \gamma_{0}-1\right)  }{\gamma_{0}}-\left(
\gamma_{\varepsilon}-1\right)  =\varepsilon.
\]
Thus, we have to verify that%
\begin{equation}
\gamma_{\varepsilon}u_{\delta,\varepsilon}^{\prime}\left(  R_{0}\right)
\leq\gamma_{0}u_{\delta,\varepsilon}^{\varepsilon}\left(  R_{0}\right)
w^{\prime}\left(  R_{0}\right)  . \label{pufa}%
\end{equation}
But%
\[
w^{\prime}\left(  R_{0}\right)  =\frac{\gamma_{\varepsilon}}{\gamma_{0}}%
\frac{1}{R_{0}^{N-1}}\int_{R_{0}}^{R}a\left(  y\right)  y^{N-1}\,dy,
\]
and so (\ref{pufa}) follows immediately from (\ref{ad}).

Taking into account (\ref{laaa}), (\ref{ann}) and (\ref{fin}),
%\marginpar{\textit{modified}}
the proof can now be ended as the proof of Theorem \ref{cc}. $\blacksquare$

\begin{remark}
Note that if we take $\delta=0$
%\marginpar{\textit{modified again}}
in the proof of Theorem \ref{rad2}, then the subsolution
vanishes at the origin. This is why we have to choose $\varepsilon>0$ and we
cannot pick $\varepsilon=0$.
\end{remark}

\begin{remark} 
Although Theorems \ref{cc} and \ref{rad2} hold in particular for $N=1$, in
this case one can obtain similar results without assuming that
$a$ is even. More precisely, if $\Omega:=\left(  \alpha,\beta\right)  $ and $\mu
\in\Omega$, a quick look at the proofs of the aforementioned theorems shows
that one can replace $B_{R_{0}}$ and $A_{R_{0},R}$ by
$\left(  \alpha,\mu\right)  $ and $\left(  \mu,\beta\right)$ respectively, in
order to reach a similar conclusion.
\end{remark}

\bigskip

\section{Proof of Theorem \ref{ti} and some corollaries}

\medskip

\textbf{Proof of Theorem \ref{ti} (i)}: First we note, as a consequence of
Theorem \ref{c1}, that $\mathcal{I}\neq\emptyset$ since $(H_{0})$ holds.
%\marginpar{\textit{modified}}

\begin{enumerate}
\item[(i1)] Let $q_{0}\in\mathcal{I}$ and $u_{0}\in\mathcal{P}^{\circ}$ be a
corresponding solution of $(P_{a,q_{0}})$. Given $q\in(q_{0},1)$, define
\[
\gamma:=\frac{1-q_{0}}{1-q}>1,\quad \text{and} \quad w:=\gamma^{\frac{-1}{1-q}}u_{0}^{\gamma}.
\]
Then, a brief computation yields that%
\begin{align*}
-\Delta w &  =-\gamma^{\frac{-1}{1-q}}\gamma\left(  \left(  \gamma-1\right)
u_{0}^{\gamma-2}\left\vert \nabla u_{0}\right\vert ^{2}+u_{0}^{\gamma-1}\Delta
u_{0}\right)  \\
&  \leq\gamma^{\frac{-1}{1-q}}\gamma u_{0}^{\gamma-1}a\left(  x\right)
u_{0}^{q_{0}}\\
&  =a\left(  x\right)  w^{q}\quad a.e.\text{ in } \Omega
\end{align*}
and $\frac{\partial w}{\partial\nu}=0$ on $\partial\Omega$. In other words,
$w$ is a subsolution of $(P_{a,q})$ belonging to $\mathcal{P}^{\circ}$. So,
recalling Remark \ref{rema} (ii), 
	% \marginpar{\textit{modified}} 
we obtain a solution $u\in\mathcal{P}^{\circ}$ of
$(P_{a,q})$, and thus $q\in\mathcal{I}$. Therefore, defining $q_{i}:=\inf$
$\mathcal{I}$, 
the first assertion of (i) follows. \newline

\item[(i2)] Since $(H_{0})$ holds, one can see by a variational approach that
%\marginpar{\textit{modified}}
$(P_{a,q})$ has a nontrivial nonnegative solution for any $0<q<1$ (see e.g.
the proof of \cite[Corollary 1.8]{krqu}), and thus $\mathcal{A}\subseteq
\mathcal{I}$.
%\marginpar{modified}
Assume now $(H_{1}^{\prime})$ and $(H_{+})$.
%\marginpar{\textit{modified}}
Let $q\in\left(  0,1\right)  $, and suppose by contradiction that there exist
$u$ and $v$ nontrivial nonnegative solutions of $(P_{a,q})$ with
$u\not \equiv v$. We claim that $u\not \equiv 0$ in $\Omega_{+}$. Indeed, if
not, then $\Delta u\geq0$ in $\Omega$ and $\frac{\partial u}{\partial\nu}=0$
on $\partial\Omega$, and therefore the maximum principle says that $u\equiv0$
in $\Omega$, which is not possible. Now, 
	% \marginpar{\textit{modified}} 
taking into account that $\Omega_{+}$ is connected (by $(H_{1}^{\prime})$), arguing as in Lemma 2.2 in
\cite{BPT} we infer that $u>0$ in 
$\Omega_{+}$. But, since the same reasoning
applies to $v$, this contradicts the uniqueness result of \cite[Theorem
3.1]{BPT} (see Remark \ref{h}). Finally, recalling that $\mathcal{A}%
\subseteq\mathcal{I}$, we deduce that $\mathcal{A}=\mathcal{I}$.
$\blacksquare$ \newline
\end{enumerate}

\medskip

\textbf{Proof of Theorem \ref{ti} (ii)}: After a dilation and a translation,
we can assume that $\Omega:=\left(  -2,2\right)  $. For any $q\in(0,1)$, we
shall construct $a\in\mathcal{C}(\overline{\Omega})$ such that $(P_{a,q})$ has
\textit{one} solution in $\mathcal{P}^{\circ}$ and \textit{two} nontrivial 
nonnegative solutions 
	%\marginpar{\textit{modified}} 
having nonempty dead cores.
This result will be proved in two parts, in accordance with the value of $q$.

\begin{enumerate}
\item First we consider $q\in\left[  \frac{1}{3},1\right)  $. We define
\[
r:=\frac{2}{1-q}\in\left[  3,\infty\right)  \quad\text{and}\quad f\left(
x\right)  :=\frac{\left(  x+1\right)  ^{r}}{r}.
\]
Note that $rq=r-2$. Let $p$ be the polynomial given by
\[
p\left(  x\right)  :=\alpha x^{3}+\beta x^{2}+\gamma x+\delta,
\]
where
\begin{gather*}
\alpha:=-\frac{2^{r-2}\left(  r+1\right)  }{3},\quad\beta:=2^{r-3}\left(
3r+1\right)  ,\\
\gamma:=-2^{r-1}\left(  r-1\right)  ,\quad\delta:=\frac{2^{r-3}}{3}(\frac
{24}{r}+5r-13).
\end{gather*}

One can verify that
\begin{equation}
p\left(  1\right)  =f\left(  1\right)  ,\quad p^{\prime}\left(  1\right)
=f^{\prime}\left(  1\right)  ,\quad p^{\prime\prime}\left(  1\right)
=f^{\prime\prime}\left(  1\right)  ,\quad p^{\prime}\left(  2\right)  =0.
\label{mm}%
\end{equation}
Moreover, it also holds that $p$ is increasing in $\left(  1,2\right)  $, and
in particular it follows that $p>0$ in $\left[  1,2\right]  $. Set%
\[
a\left(  x\right)  :=%
\begin{cases}
-\left(  r-1\right)  r^{q} & \text{if }x\in\left[  0,1\right]  ,\\
-\frac{p^{\prime\prime}\left(  x\right)  }{\left[  p\left(  x\right)  \right]
^{q}} & \text{if }x\in\left[  1,2\right]  ,\\
a\left(  -x\right)  & \text{if }x\in\left[  -2,0\right]  ,
\end{cases}
\]
and observe that $a\in\mathcal{C}(\overline{\Omega})$ since (recall that
$rq=r-2$)%
\[
-\frac{p^{\prime\prime}\left(  1\right)  }{\left[  p\left(  1\right)  \right]
^{q}}=-\frac{f^{\prime\prime}\left(  1\right)  }{\left[  f\left(  1\right)
\right]  ^{q}}=-\left(  r-1\right)  r^{q}.
\]
Also, since $p>0$ in $\left[  1,2\right]  $, it follows from the definition
that $a$ changes sign in $\left(  1,2\right)  $. Furthermore,%
\begin{align*}
\int_{1}^{2}a  &  =-\int_{1}^{2}\frac{p^{\prime\prime}\left(  x\right)
}{\left[  p\left(  x\right)  \right]  ^{q}}=-\left[  \left.  \frac{p^{\prime
}\left(  x\right)  }{\left[  p\left(  x\right)  \right]  ^{q}}\right\vert
_{1}^{2}+q\int_{1}^{2}\frac{\left[  p^{\prime}\left(  x\right)  \right]  ^{2}%
}{\left[  p\left(  x\right)  \right]  ^{q+1}}\right] \\
&  <\left.  \frac{p^{\prime}\left(  x\right)  }{\left(  \left[  p\left(
x\right)  \right]  \right)  ^{q}}\right\vert _{2}^{1}=\frac{p^{\prime}\left(
1\right)  }{\left[  p\left(  1\right)  \right]  ^{q}}=\frac{f^{\prime}\left(
1\right)  }{\left[  f\left(  1\right)  \right]  ^{q}}=2r^{q}%
\end{align*}
and hence%
\[
\int_{0}^{2}a<2r^{q}-\left(  r-1\right)  r^{q}\leq0
\]
%\marginpar{\textit{modified}}
since $r\geq3$. Therefore, $\int_{\Omega}a<0$. Define now
\[
u_{1}\left(  x\right)  :=%
\begin{cases}
0 & \text{if }x\in\left[  -2,-1\right]  ,\\
f\left(  x\right)  & \text{if }x\in\left[  -1,1\right]  ,\\
p\left(  x\right)  & \text{if }x\in\left[  1,2\right]  ,
\end{cases}
\]
and $u_{2}\left(  x\right)  :=u_{1}\left(  -x\right)  $. Taking into account
(\ref{mm}), we see that $u_{1},u_{2}\in\mathcal{C}^{2}(\overline{\Omega})$.
Moreover, one can see that 
	% \marginpar{\textit{modified}} 
$u_{1}$ and $u_{2}$ are two distinct nonnegative nontrivial
solutions of the problem 
\begin{equation}
\left\{
\begin{array}
[c]{lll}%
-u^{\prime\prime}=a(x)u^{q} & \mathrm{in} & \Omega,\\
u^{\prime}=0 & \mathrm{on} & \partial\Omega.
\end{array}
\right.  \label{eje}%
\end{equation}

Now, a simple integration by parts shows that
\[
\max\left(  u_{1}\left(  x\right)  ,u_{2}\left(  x\right)  \right)  =\left\{
\begin{array}
[c]{lll}%
u_{2}(x) & \mathrm{in} & \left[  -2,0\right]  ,\\
u_{1}(x) & \mathrm{in} & \left[  0,2\right]  ,
\end{array}
\right.
\]
is a strictly positive weak subsolution of (\ref{eje}). Thus, since
$\int_{\Omega}a<0$, by Remark \ref{rema} (ii) there exist arbitrary large
supersolutions of (\ref{eje}) and we then obtain a solution $u\in
\mathcal{P}^{\circ}$ of (\ref{eje}). It follows that $q\in\mathcal{I}$, but
$q\not \in \mathcal{A}$, since $u_{1}$ and $u_{2}$ are nontrivial nonnegative
solutions having 
nonempty dead cores. \newline

\item Now we consider $q\in\left(  0,\frac{1}{3}\right)  $. We proceed as
above, the only difference being the definition of $p$. For $K>0$, let%
\[
p\left(  x\right)  =p_{K}\left(  x\right)  :=\alpha x^{4}+\beta x^{3}+\gamma
x^{2}+\delta x+\mu,
\]
where
\begin{gather*}
\alpha:=3(\frac{2^{r}}{r}-K)+2^{r-3}\left(  r+7\right)  ,\quad\beta
:=-16(\frac{2^{r}}{r}-K)+2^{r-3}\left(  -6r-38\right)  ,\\
\gamma:=30(\frac{2^{r}}{r}-K)+2^{r-3}\left(  13r+71\right)  ,\quad
\delta:=-24(\frac{2^{r}}{r}-K)-2^{r-3}\left(  12r+52\right)  ,\\
\mu:=8\frac{2^{r}}{r}-7K+2^{r-3}\left(  4r+12\right)  .
\end{gather*}
One can check that
\begin{gather}
p\left(  1\right)  =f\left(  1\right)  ,\quad p^{\prime}\left(  1\right)
=f^{\prime}\left(  1\right)  ,\quad p^{\prime\prime}\left(  1\right)
=f^{\prime\prime}\left(  1\right)  ,\label{asz}\\
p^{\prime}\left(  2\right)  =0,\quad p\left(  2\right)  =K.\nonumber
\end{gather}
We observe that $p>0$ in $\left[  1,2\right]  $. Indeed, since
\[
p\left(  1\right)  ,p^{\prime}\left(  1\right)  ,p^{\prime\prime}\left(
1\right)  ,p\left(  2\right)  >0=p^{\prime}\left(  2\right)  ,
\]
if $p\leq0$ somewhere, then $p^{\prime\prime}$ would vanish at least at three
points in $\left(  1,2\right)  $, which is impossible since $p^{\prime\prime}$
has degree $2$. It follows that
\[
a_{K}\left(  x\right)  :=-\frac{p_{K}^{\prime\prime}\left(  x\right)
}{\left[  p_{K}\left(  x\right)  \right]  ^{q}}\in\mathcal{C}(\left[
1,2\right]  ).
\]
%\marginpar{\textit{modified, from here, the rest of the proof}}
We claim now that for all $K>0$ large enough, $a_{k}$ changes sign in $\left(
1,2\right)  $ and $\int_{1}^{2}a_{K}<0$. Indeed, $a_{K}(1)<0$, and for $K$
sufficiently large we have $p_{K}^{\prime\prime}(2)<0$, so that $a_{K}(2)>0$.
Hence, the first assertion follows. To show the second one, we first note that%
\begin{equation}
\int_{1}^{2}a_{K}=-K^{1-q}\int_{1}^{2}\frac{p_{K}^{\prime\prime}\left(
x\right)  }{K}\left(  \frac{K}{p_{K}\left(  x\right)  }\right)  ^{q}%
,\label{limo}%
\end{equation}%
\[
\lim_{K\rightarrow\infty}\frac{p_{K}\left(  x\right)  }{K}=-3x^{4}%
+16x^{3}-30x^{2}+24x-7=\left(  7-3x\right)  \left(  x-1\right)  ^{3}:=g\left(
x\right)  ,
\]
and
\[
\lim_{K\rightarrow\infty}\frac{p_{K}^{\prime\prime}\left(  x\right)  }%
{K}=-36x^{2}+96x-60=12\left(  5-3x\right)  \left(  x-1\right)  :=12h\left(
x\right)  .
\]
Define%
\[
H(x):=-x^{3}+4x^{2}-5x+2=(2-x)(x-1)^{2}>0\quad\text{in }(1,2).
\]
Then, 
	% \marginpar{\textit{modified}} 
$H^{\prime}=h$. Also, since $2-3q>0$ we see that
$\displaystyle\lim_{x\rightarrow1^{+}}H(x)g(x)^{-q}=0$. Therefore, an
integration by parts yields that
\begin{align}
\int_{1}^{2}\frac{h\left(  x\right)  }{g^{q}\left(  x\right)  } &
=H(x)g(x)^{-q}\big{\vert}_{1}^{2}-\int_{1}^{2}H(x)\left(  g(x)^{-q}\right)
^{\prime}\label{limon}\\
&  =q\int_{1}^{2}H(x)g(x)^{-\left(  q+1\right)  }g^{\prime}(x)>0\nonumber
\end{align}
because $g^{\prime}>0$ in $(1,2)$. It follows from (\ref{limo}) and
(\ref{limon}) that
\[
\lim_{K\rightarrow\infty}\int_{1}^{2}a_{K}=-\infty,
\]
and therefore the claim is proved. We can then fix some $K>0$ such that
$a_{K}$ changes sign in $\Omega$ and $\int_{\Omega}a_{K}<0$, and thus the
proof can be completed as in the previous case. $\blacksquare$
\end{enumerate}

\medskip

\begin{remark}
Let us point out that, by the uniqueness results in \cite{BPT}, for every
$q\in\left(  0,1\right)  $, the problem $(P_{a,q})$ with the weight $a_{q}$
constructed in the above proof has \textit{exactly} three (nontrivial)
nonnegative solutions. Indeed, one can verify that $\Omega_{+}(a_{q})$ has
exactly two connected components (taking $K$ large if $q<1/3$), say
$\mathcal{O}_{1}$ and $\mathcal{O}_{2}$. Now, by \cite[Theorem 3.1]{BPT} there
exists at most one nonnegative solution which is positive in $\mathcal{O}_{1}$
and zero in $\mathcal{O}_{2}$, and vice-versa. Also, by the aforementioned
theorem, there exists at most one nonnegative solution which is positive in
both $\mathcal{O}_{1}$ and $\mathcal{O}_{2}$. Since the nontrivial nonnegative
solutions $u$ satisfy that either $u>0$ in $\mathcal{O}_{i}$ or $u\equiv0$ in
$\mathcal{O}_{i}$ (see \cite[Lemma 2.2]{BPT}), our assertion follows because
from the maximum principle we deduce that there is no nontrivial nonnegative
solution vanishing in both $\mathcal{O}_{1}$ and $\mathcal{O}_{2}$.\newline
Let us also remark that the solution in $\mathcal{P}^{\circ}$ is even: indeed,
if not, we would have four nontrivial nonnegative solutions. Summing up, for
%\marginpar{\textit{modified}}
this family of even weights $a_{q}$, there exist two (nontrivial) noneven
nonnegative solutions with nonempty dead cores, and one even solution in
$\mathcal{P}^{\circ}$.
\end{remark}

\begin{remark}
Let $a_{q}$ be as in the first case of the proof of Theorem \ref{ti} (ii), but now with $q\in\left[  0,1\right)  $. A quick look at the aforementioned proof
shows that $\int_{\Omega}a_{q}>0$ for $q>0$ close enough to $0$. Indeed, this
follows easily from the fact that $\int_{\Omega}a_{0}=1$. Furthermore, for
such $q$'s, reasoning as therein we obtain two (nontrivial) nonnegative
solutions of $(P_{a_q,q})$. In other words, this result shows that, unlike for the
existence of positive solutions, the condition $\int_{\Omega}a<0$ is
\textit{not necessary} in order to have existence of (nontrivial) nonnegative
solutions of $(P_{a,q})$. Let us add that this matter has already been noted in
\cite[Section 2.3]{BPT}.
\end{remark}

\medskip

As an immediate consequence of Theorems \ref{rad2}
%\marginpar{\textit{modified}}
and \ref{ti} (i), we have the following result:

\begin{corollary}
\label{cor:KN}
%\marginpar{\textit{modified}}
Let $\Omega:=B_{R}$ and $a\in L^{\infty}\left(  \Omega\right)  $\textit{ be a
radial function such that }$\int_{\Omega}a<0$ and $a\geq0$\textit{\ in
}$A_{R_{0},R}$ for some $R_{0}\in\left(  0,R\right)  $. Then, 
	% \marginpar{\textit{modified below}} 
\[
\left(  \frac{1-KN}{1-KN+2K},1\right)  \subseteq\mathcal{I},\quad
\text{where\quad}K=K(a):=\frac{\int_{A_{R_{0},R}}a}{\omega_{N-1}R_{0}%
^{N}\left\Vert a^{-}\right\Vert _{L^{\infty}(B_{R_{0}})}}.
\]
Moreover,
\[
\left(  \frac{1-KN}{1-KN+2K},1\right)  \subseteq\mathcal{A}%
\]
if $a\leq0$ in $B_{R_{0}}$.
\end{corollary}

\textbf{Proof}. Since $\int_{\Omega}a<0$, a direct computation gives that
$KN<1$. Let $q\in\left(  \frac{1-KN}{1-KN+2K},1\right)  $. Then one can check
that (\ref{sipi}) is satisfied and thus there exists $u\in\mathcal{P}^{\circ}$
solution of $(P_{a,q})$, so that $q\in\mathcal{I}$. The last assertion of the
corollary is now immediate from Theorem \ref{ti} (i2). $\blacksquare$\newline

\begin{remark}
\label{rem:I01}

Let us point out 
	% \marginpar{\textit{modified, and also the next line}} 
that $\mathcal{I}(a)$ may approach the whole interval $(0,1)$ as
the coefficient $a$ varies. To show this, we may use either the sub and
supersolutions method (Corollary \ref{cor:KN}), or a bifurcation analysis
(Proposition \ref{prop:q0}). Let us also add that, however, we believe that
there is no $a$ such that $\mathcal{I}\left(  a\right)  =\left(  0,1\right)
$, but we are not able to prove it.\strut

\begin{enumerate}
\item[(i)] Given any fixed $q_{0}\in\left(  0,1\right)  $, Corollary
\ref{cor:KN} provides some cases in which $\left(  q_{0},1\right)
\subseteq\mathcal{I}=\mathcal{A}$. Indeed, in order to see this it suffices to
find $a$ such that $K(a)$ satisfies $1>K(a)N\approx1$. One may take for
instance $\Omega:=B_{1}$ and
\[
a\left(  x\right)  :=\sigma\chi_{A_{\frac{1}{2},1}}\left(  x\right)
-\chi_{B_{\frac{1}{2}}}\left(  x\right)  ,\quad\text{ for }x\in\Omega,
\]
%\marginpar{\textit{modified}}
where $\sigma>0$. Since $K(a)N=\sigma\left(  2^{N}-1\right)  $, it is easy to
choose $\sigma$ adequately.

\item[(ii)] Let $a\in L^{\infty}(\Omega)$ be given by $a(x):=\sigma
\chi_{\Omega_{1}}(x)-\chi_{\Omega_{2}}(x)a_{2}(x)$,
%\marginpar{\textit{modified}}
where $\sigma>0$, $a_{2}\geq0$, and $\Omega_{1},\Omega_{2}$ are disjoint
subsets of $\Omega$ such that $\int_{\Omega}a=0$. Then, for any $\varepsilon
>0$ small, we see that $a-\varepsilon$ changes sign, $\int_{\Omega
}(a-\varepsilon)<0$, and $\Omega_{+}(a-\varepsilon)=\Omega_{+}(a)$. By
combining
%Theorem 1.9 in\cite{krqu},
Theorem \ref{ti} (i) and Proposition \ref{prop:q0}, we see that $\mathcal{I}%
_{a-\varepsilon}$ approaches $(0,1)$ as $\varepsilon\rightarrow0^{+}$.
Additionally, if $\Omega_{+}(a)$ satisfies $(H_{1}^{\prime})$ and $(H_{+})$,
then $\mathcal{A}_{a-\varepsilon}=\mathcal{I}_{a-\varepsilon}$ approaches
$(0,1)$ as $\varepsilon\rightarrow0^{+}$. \newline
\end{enumerate}
\end{remark}

%Combining Theorem \ref{ta}, Theorem \ref{ti}(1), and Proposition \ref{prop:q0} provides the following result on the properties of the sets $\mathcal{A}=\mathcal{A}_a$ and $\mathcal{I}=\mathcal{I}_{a}$, which corresponds to Remark \ref{rem:I01}.

%\begin{corollary}
%Assume that $(H_1)$ holds. Let $b \in L^r(\Omega)$, $r>N$, be such that $b\not\equiv 0$, and $\int_\Omega b = 0$. Then, we have $\mathcal{I}_{b-\varepsilon} \rightarrow (0,1)$ as $\varepsilon \rightarrow 0^+$. Additionally if $b \in L^\infty (\Omega)$, then we have $\mathcal{A}_{b-\varepsilon} = \mathcal{I}_{b-\varepsilon} \rightarrow (0,1)$ as $\varepsilon \rightarrow 0^+$, whenever $\Omega_+$ is conneceted, and $(H_+)$ holds. \\
%\end{corollary}

\section{Proof of Theorem \ref{dc}}

\label{sec:core} \textbf{Proof of Theorem \ref{dc} (i):} Let $\Omega_{0}$ be a
tubular neighborhood of $\partial\Omega$ such that $a>0$ $a.e.$ in $\Omega
_{0}$, with smooth boundary $\partial\Omega_{0}=\partial\Omega\cup\Gamma_{0}$,
where $\Gamma_{0}=\partial\Omega_{0}\cap\Omega$. We consider the following
concave mixed problem
\begin{equation}%
\begin{cases}
-\Delta v=a^{+}(x)v^{q} & \mbox{ in }\Omega_{0},\\
\frac{\partial v}{\partial\nu}=0 & \mbox{ on }\partial\Omega,\\
v=0 & \mbox{ on }\Gamma_{0}.
\end{cases}
\label{mixpp}%
\end{equation}
Proceeding as in \cite[Lemma 3.3]{ABC94}, we see that the comparison principle
holds for \eqref{mixpp}, i.e., $\underline{v}\leq\overline{v}$ on
$\overline{\Omega_{0}}$ for any nonnegative supersolution $\overline{v}$ and
subsolution $\underline{v}$ of \eqref{mixpp} such that $\overline
{v},\underline{v}>0$ in $\Omega_{0}$.

Let $u$ be a nontrivial nonnegative solution of $(P_{a,q})$. Then $u$ is a
supersolution of \eqref{mixpp}. In addition, $u>0$ in $\Omega_{0}$.
%\marginpar{\textit{modified}}
Indeed, recalling $(H_{1}^{\prime})$, we observe that
\[
0<\int_{\Omega}|\nabla u|^{2}=\int_{\Omega}a(x)u^{q+1}\leq\int_{\text{supp
}a^{+}}a(x)u^{q+1}=\int_{\Omega_{+}}a(x)u^{q+1}.
\]
It follows that $u\not \equiv 0$ in $\Omega_{+}$.
%\marginpar{\textit{modified}%}
Also, since $\Omega_{+}$ is connected (by $(H_{1}^{\prime})$), the strong
maximum principle yields that $u>0$ in $\Omega_{+}$, and consequently $u>0$ in
$\Omega_{0}$ as claimed.

On the other hand, in order to construct a subsolution, we consider the
following mixed eigenvalue problem:
\begin{equation}%
\begin{cases}
-\Delta\psi=\sigma a^{+}(x)\psi & \mbox{ in }\Omega_{0},\\
\frac{\partial\psi}{\partial\nu}=0 & \mbox{ on }\partial\Omega,\\
\psi=0 & \mbox{ on }\Gamma_{0}.
\end{cases}
\label{mixep}%
\end{equation}
By $\sigma_{1}>0$, we denote the smallest eigenvalue of \eqref{mixep}, and by
$\psi_{1}$ an eigenfunction associated to $\sigma_{1}$ satisfying $\psi_{1} >
0$ on $\overline{\Omega_{0}} \setminus\Gamma_{0}$. Then, we see that
$\varepsilon_{q}\psi_{1}$ is a subsolution of \eqref{mixep} for some
$\varepsilon_{q}>0$ small. By the comparison principle, we deduce that
$\varepsilon_{q}\psi_{1}\leq u$ on $\overline{\Omega_{0}}$, from which the
desired conclusion follows. $\blacksquare$\newline

The following result will be used in the proof of Theorem \ref{dc} (ii):
%\marginpar{\textit{Lemma 5.1 is now outside the proof of Theorem \ref{dc}(ii)}}

\begin{lemma}
\label{lem:bound} Under the conditions of Theorem \ref{dc} (ii), let
$\overline{q}\in(0,1)$ and $\delta_{2} > 0$ such that $\int_{\Omega}%
a_{\delta_{2}}< 0$. Then, there exists 
	% \marginpar{\textit{modified}}  
$C>1$ such that $\| u\|_{L^{\infty
}(\Omega)} < C$ for all nonnegative solutions $u$ of $(P_{a_{\delta},q})$ with
$q \in(0, \overline{q}]$ and $\delta\geq\delta_{2}$. \newline
\end{lemma}

\textbf{Proof. }Let $u$ be a nonnegative solution of $(P_{a_{\delta},q})$ with
$q\in(0,\overline{q}]$ and $\delta\geq\delta_{2}$. Then $u$ is a subsolution
of $(P_{a_{\delta},q})$ for $\delta=\delta_{2}$.
%\[
%\left\{
%\begin{array}
%[c]{lll}%
%-\Delta v=(b_{1}(x)-\delta_{2}b_{2}(x))v^{q} & \mathrm{in} & \Omega,\\
%\frac{\partial v}{\partial\nu}=0 & \mathrm{on} & \partial\Omega.
%\end{array}
%\right.  % \label{pdel0}%
%\]
%\marginpar{\textit{modified, this is why we had to modify }$H_{2}$}
In view of Remark \ref{rema} (ii), we can construct a supersolution $w$ of
$(P_{a_{\delta_{2}},q})$ such that $u\leq w$. Hence, the sub and
supersolutions method ensures the existence of a nonnegative solution $v$ of
$(P_{a_{\delta_{2}},q})$ such that $u\leq v\leq w$. By Proposition \ref{t2.1}
(ii), $(P_{a_{\delta_{2}},q})$ has an \textit{a priori} bound for nonnegative
solutions in $L^{\infty}(\Omega)$, which is uniform in $q\in(0,\overline{q}]$.
The lemma 
	% \marginpar{\textit{modified}} 
now follows. $\blacksquare$\newline

\noindent\textbf{Proof of Theorem \ref{dc} (ii): }
%We first establish an \textit{a priori} bound in $L^{\infty}(\Omega)$ for nonnegative solutions of $(P_{a_{\delta},q})$. This bound shall be uniform in $q \in(0, \overline{q}]$ and $\delta\to\infty$. Now, W
We proceed as in the proofs of \cite[Theorem 1(iv)]{GMRS07} and \cite[Theorem
3.1]{FP84}.
%to prove Theorem \ref{dc} (ii).
Let $\overline{q}\in(0,1)$ and $\delta_{2}>0$ be the constant given by Lemma
\ref{lem:bound}, and $u$ be a nontrivial nonnegative solution of
$(P_{a_{\delta},q})$ with $q\in(0,\overline{q}]$ and $\delta\geq\delta_{2}$.
Given $\sigma>0$, we pick $a_{0}>0$ such that
\begin{equation}
b_{2}(z)\geq a_{0}\quad\mbox{ for all }z\in G_{\sigma/2},\label{b2a0}%
\end{equation}
where we have used the continuity of $b_{2}$. Let us fix $x\in
G_{\sigma/2}$, and consider $d(x):=\text{dist}(x,\partial G_{\sigma/2})$,
where $d(x)>0$ since $G_{\sigma/2}$ is open. We then define
\begin{equation}
v_{1}(y):=d(x)^{-\alpha}u(x+d(x)y),\quad\text{for }\left\vert y\right\vert
\leq1.\label{def:v1}%
\end{equation}
Let $\alpha:=2/(1-q)$, so that $2-\alpha+\alpha q=0$. If $\left\vert
y\right\vert <1$ then, using \eqref{b2a0}, a brief computation yields
\begin{align*}
-\Delta v_{1}(y) &  =d(x)^{2-\alpha}a_{\delta}(x+d(x)y)\,u(x+d(x)y)^{q}\\
&  =-\delta b_{2}(x+d(x)y)\,v_{1}(y)^{q}\\
&  \leq-\delta a_{0}\,v_{1}(y)^{q}.
\end{align*}
Here, we have used the fact that $x+d(x)y\in G_{\sigma/2}$. If $\left\vert
y\right\vert =1$, then we have
\begin{equation}
v_{1}(y)\leq d(x)^{-\alpha}C,\label{c5}%
\end{equation}
where $C>1$ is provided by Lemma \ref{lem:bound}.

Given $\varepsilon>0$, we now consider the problem
\[
\left\{
\begin{array}
[c]{lll}%
-\Delta v=-\delta a_{0}v^{q} & \mathrm{in} & B_{1},\\
v=\varepsilon & \mathrm{on} & \partial B_{1}.
\end{array}
\right.  \leqno{(Q_{\delta, \varepsilon})}
\]
We observe from \eqref{c5} that $v_{1}$ is a subsolution of $(Q_{\delta
,\varepsilon})$ if
\begin{equation}
d(x)^{-\alpha}C\leq\varepsilon. \label{C1}%
\end{equation}

Next, we construct a supersolution of $(Q_{\delta,\varepsilon})$. For
$r=\left\vert y\right\vert $, we define
\begin{equation}
z_{1}(r):=\left\{
\begin{array}
[c]{ll}%
0, & 0\leq r\leq\frac{1}{2},\\
A(r-\frac{1}{2})^{\alpha}, & \frac{1}{2}<r\leq1,
\end{array}
\right.  \label{def:z1}%
\end{equation}
where $A$ is a positive constant to be determined. Since $\alpha>2$, we have
$z_{1}\in\mathcal{C}^{2}(\overline{B_{1}})$, and in addition,
\begin{align*}
\Delta z_{1}=z_{1}^{\prime\prime}+\frac{N-1}{r}z_{1}^{\prime}  &
=A\alpha(\alpha-1)(r-\frac{1}{2})^{\alpha-2}+\frac{N-1}{r}A\alpha(r-\frac
{1}{2})^{\alpha-1}\\
&  \leq A\alpha(\alpha-1)(r-\frac{1}{2})^{\alpha-2}+(N-1)A\alpha(r-\frac{1}%
{2})^{\alpha-2}\\
&  \leq\delta a_{0}\left(  A(r-\frac{1}{2})^{\alpha}\right)  ^{q}=\delta
a_{0}z_{1}^{q}\quad\text{ for }\frac{1}{2}<r<1,
\end{align*}
if $A\alpha(\alpha-1)+(N-1)A\alpha\leq\delta a_{0}A^{q}$, i.e.,
\begin{equation}
A\leq\left(  \frac{\delta a_{0}}{\alpha(\alpha-1)+(N-1)\alpha}\right)
^{\frac{1}{1-q}}. \label{c2}%
\end{equation}
Moreover, we note that
\[
z_{1}(1)=A\left(  \frac{1}{2}\right)  ^{\alpha}\geq\varepsilon,
\]
provided that
\begin{equation}
A\geq2^{\alpha}\varepsilon. \label{c3}%
\end{equation}
Hence, from \eqref{C1}, \eqref{c2} and \eqref{c3}, it follows that if
\[
2^{\alpha}d(x)^{-\alpha}C\leq2^{\alpha}\varepsilon\leq A\leq\left(
\frac{\delta a_{0}}{\alpha(\alpha-1)+(N-1)\alpha}\right)  ^{\frac{1}{1-q}},
\]
i.e.,
\[
d(x)\geq2C^{\frac{1}{\alpha}}\left(  \frac{\alpha(\alpha-1)+(N-1)\alpha
}{\delta a_{0}}\right)  ^{\frac{1}{2}},
\]
then $v_{1}$ is a subsolution of $(Q_{\delta,\varepsilon})$, and in addition,
$z_{1}$ is a supersolution of $(Q_{\delta,\varepsilon})$. Since $2<\alpha
\leq2/(1-\overline{q})=:\overline{\alpha}$, this occurs for some
$\varepsilon=\varepsilon_{\delta,x}$ if
\begin{equation}
d(x)\geq2\left(  \frac{\overline{\alpha}(\overline{\alpha}-1)+(N-1)\overline
{\alpha}}{C^{-1}\delta a_{0}}\right)  ^{\frac{1}{2}}=:d_{\delta}. \label{c4}%
\end{equation}

Now, using the comparison principle for $(Q_{\delta,\varepsilon})$
%\marginpar{\textit{modified}}
(which is deduced from the weak maximum principle) we derive that $v_{1}\leq
z_{1}$, so that
\[
d(x)^{-\alpha}u(x)=v_{1}(0)\leq z_{1}(0)=0,
\]
and consequently, $u(x)=0$. Therefore, we have proved that if $x\in
G_{\sigma/2}$ satisfies \eqref{c4}, then $u(x)=0$ for any nontrivial
nonnegative solution $u$ of $(P_{a_{\delta},q})$ with $q\in(0,\overline{q}]$.
Since $d_{\delta}$ in \eqref{c4} does not depend on $u$ or $q$, and converges
to $0$ as $\delta\rightarrow\infty$, we have the desired conclusion.
$\blacksquare$ \newline

%%%

\subsection*{Acknowledgements}

U. Kaufmann was partially supported by Secyt-UNC. H. Ramos Quoirin was
supported by Fondecyt 1161635. K. Umezu was supported by JSPS KAKENHI Grant
Number 15K04945.

%%%

\end{document}